\documentclass[amscd,amssymb,11pt]{amsart}

\usepackage[all]{xy}

\newtheorem{thm}{Theorem}[section]
\newtheorem{lem}[thm]{Lemma}
\newtheorem{cor}[thm]{Corollary}
\newtheorem{prop}[thm]{Proposition}

\theoremstyle{definition}

\newtheorem{defn}[thm]{Definition}
\newtheorem{defns}[thm]{Definitions}

\newtheorem{ex}[thm]{Example}

\theoremstyle{remark}

\newtheorem{rem}[thm]{Remark}
\newtheorem{rems}[thm]{Remarks}

\numberwithin{equation}{section}

\newcommand{\thmref}[1]{Theorem~\ref{#1}}
\newcommand{\corref}[1]{Corollary~\ref{#1}}
\newcommand{\secref}[1]{\S\ref{#1}}

\newcommand{\propref}[1]{Proposition~\ref{#1}}
\newcommand{\lemref}[1]{Lemma~\ref{#1}}
\newcommand{\appref}[1]{Appendix~\ref{#1}}
\newcommand{\exref}[1]{Example~\ref{#1}}

\newcommand{\hocolim}{\operatorname*{hocolim}}
\newcommand{\holim}{\operatorname*{holim}}
\newcommand{\colim}{\operatorname*{colim}}
\newcommand{\totfib}{\operatorname*{TotFib}}

\newcommand{\Tor}{\operatorname{Tor}}

\newcommand{\A}{{\mathcal  A}}
\newcommand{\B}{{\mathcal  B}}
\newcommand{\C}{{\mathcal  C}}
\newcommand{\Sp}{{\mathcal  S}}

\newcommand{\Alg}{{\mathcal Alg}}

\newcommand{\Algg}{{\mathcal Alg}^{\prime}}

\newcommand{\Module}{{\mathcal Mod}}

\newcommand{\Z}{{\mathbb  Z}}

\newcommand{\Q}{{\mathbb  Q}}
\newcommand{\PP}{{\mathbb P}}

\newcommand{\Sinfty}{\Sigma^{\infty}}
\newcommand{\Oinfty}{\Omega^{\infty}}

\newcommand{\sm}{\wedge}

\newcommand{\ra}{\rightarrow}
\newcommand{\xra}{\xrightarrow}
\newcommand{\la}{\leftarrow}
\newcommand{\xla}{\xleftarrow}

\newcommand{\hra}{\hookrightarrow}

\begin{document}

\title[Periodic homology of infinite loopspaces]{Localization of Andr\'{e}--Quillen--Goodwillie towers, and the periodic homology of infinite loopspaces}             

\author[Kuhn]{Nicholas J.~Kuhn}                                  
\address{Department of Mathematics \\ University of Virginia \\ Charlottesville, VA 22903}    
\email{njk4x@virginia.edu}
\thanks{This research was partially supported by a grant from the National Science Foundation}     
 
\date{June 3, 2003.}


\begin{abstract}   
Let $K(n)$ be the $n^{th}$ Morava K--theory at a prime $p$, and let $T(n)$ be the telescope of a $v_n$--self map of a finite complex of type $n$.  In this paper we study the $K(n)_*$--homology of  $\Omega^{\infty} X$, the $0^{th}$ space of a spectrum $X$, and many related matters.  

We give a sampling of our results.

Let $\PP X$ be the free commutative $S$--algebra generated by $X$: it is weakly equivalent to the wedge of all the extended powers of $X$.  We construct a natural map 
$$ s_n(X): L_{T(n)}\PP(X) \ra L_{T(n)}\Sinfty(\Oinfty X)_+$$
of commutative algebras over the localized sphere spectrum $L_{T(n)}S$.
The induced map of commutative, cocommutative $K(n)_*$--Hopf algebras
$$ s_{n}(X)_*: K(n)_*(\PP X) \ra K(n)_*(\Oinfty X),$$
satistfies the following properties.

It is always monic.

It is an isomorphism if $X$ is $n$--connected, $\pi_{n+1}(X)$ is torsion, and $T(i)_*(X) = 0$ for $1 \leq i \leq n-1$. It is an isomorphism only if $K(i)_*(X) = 0$ for $1 \leq i \leq n-1$. 

It is universal: the domain of $s_n(X)_*$ preserves $K(n)_*$--isomorphisms, and if $F$ is any functor preserving $K(n)_*$--isomorphisms, then any natural transformation $F(X) \ra K(n)_*(\Oinfty X)$ factors uniquely through $s_n(X)_*$.

The construction of our natural transformation uses the telescopic functors constructed and studied previously by Bousfield and the author, and thus depends heavily on the Nilpotence Theorem of Devanitz, Hopkins, and Smith.  Our proof that $s_n(X)_*$ is always monic uses Topological Andr\'{e}--Quillen Homology and Goodwillie Calculus in nonconnective settings.

\end{abstract}
                                                  
\maketitle

\section{Introduction and main results} \label{introduction}
 
In algebraic topology, homotopical aspects of topological spaces are studied by means of generalized homology and cohomology theories.  Such theories are themselves determined by spectra, the objects of the stable category.  
One can then pass between the worlds of unstable and stable homotopy by means of the adjoint pair of functors $(\Sinfty, \Oinfty)$, where $\Sinfty Z$ denotes the suspension spectrum of a based space $Z$, and $\Oinfty X$ denotes the $0^{th}$ infinite loopspace of the spectrum $X$.  

Though $\Sinfty$ preserves homology (and cohomology), the homological behavior of $\Oinfty$ is much more subtle,  and one has the basic problem: 
given a spectrum $E$, to what extent, and in what way, is $E_*(\Oinfty X)$ determined by $E_*(X)$?  

There is a related, more subtle, problem: to what extent, and in what way, is $L_E\Sinfty \Oinfty X$ determined by $L_EX$?  Here $L_E$ denotes Bousfield localization with respect to $E_*$.

In this paper, we develop new techniques allowing for a thorough study of these questions when $E_*$ is a periodic homology theory.  The key to our methods is to combine two of the major strands of homotopy theory of the past two decades: the flowering of powerful new techniques in homotopical algebra, many following the conceptual model offered by T. Goodwillie's calculus of functors \cite{goodwillie1, goodwillie2, goodwillie3}, and the deepening of our understanding of homotopy as organized from the chromatic point of view, in the wake of the Nilpotence Theorem of E. Devanitz, M. Hopkins, and J. Smith \cite{dhs}.

The tools from modern homotopical algebra that we use are  Topological Andr\'{e}--Quillen homology of $E_{\infty}$--ring spectra, as developed via the Goodwillie calculus framework, together with a good theory of Bousfield localization of structured objects.  These concepts require that we work within a nice model category of spectra.  Thus, for us, spectra will mean objects in $\Sp$, the category of $S$--modules as in  \cite{ekmm}, and so, e.g., commutative $S$--algebras will serve as $E_{\infty}$--ring spectra.

The input from chromatic homotopy theory comes from our use of the telescopic functors, constructed by Bousfield and the author in \cite{k2}, \cite{bousfield1}, and \cite{bousfield4}, which factor certain periodic localization functors through $\Oinfty$.  

The classic stable splitting of the space $\Oinfty \Sinfty Z$ \cite{kahn} provides a model for our main results when presented as follows.  Let $Z_+$ denote the union of a space $Z$ with a disjoint basepoint. If $X$ is an $S$--module, $\Sinfty (\Oinfty X)_+$ is naturally a commutative $S$--algebra.  Another example is $\PP X$, the free commutative algebra generated by $X$.  This is weakly equivalent to the wedge, running over $r \geq 0$, of $D_rX = E\Sigma_{r+} \sm_{\Sigma_r} X^{\sm r}$, the $r^{th}$ extended power of $X$.  Then, for all spaces $Z$, there is a natural map of commutative $S$--algebras
$$ s(Z): \PP(\Sinfty Z) \ra \Sinfty(\Oinfty \Sinfty Z)_+$$
satisfying the following properties. \\

\noindent (1) $s(Z)_*: E_*(\PP(\Sinfty Z)) \ra E_*(\Oinfty \Sinfty Z)$ is monic for all theories $E_*$. \\

\noindent (2) $s(Z)$ is an equivalence if $Z$ is connected. \\

Now let $K(n)$ be the $n^{th}$ Morava K--theory spectrum at a fixed prime $p$, and $T(n)$ its `telescopic' variant: the telescope of a $v_n$--self map of a finite complex of type $n$. (See \cite{ravenel} for background on this material.) The Telescope Conjecture is the statement that $\langle K(n) \rangle \geq \langle T(n) \rangle$, i.e. $T(n)_*$--acyclics are $K(n)_*$--acyclics; in any case, the converse holds, so that $L_{K(n)}L_{T(n)} = L_{K(n)}$.

We show that, for all spectra $X$, there is a natural map of commutative $L_{T(n)}S$--algebras
$$ s_n(X): L_{T(n)}\PP(X) \ra L_{T(n)}\Sinfty(\Oinfty X)_+$$
satisfying the following properties. \\

\noindent (1) $s_n(Z)_*: E_*(\PP(X)) \ra E_*(\Oinfty X)$ is monic for all $X$, if $\langle T(n) \rangle \geq \langle E \rangle$. \\

\noindent (2) $s_n(X)$ is an equivalence if $X$ is suitably connected and $T(i)_*(X) = 0$ for $1 \leq i \leq n-1$. It is an equivalence only if $K(i)_*(X) = 0$ for $1 \leq i \leq n-1$. \\

\noindent (3) $s_n$ is universal in the sense that any natural transformation from a functor invariant under $T(n)_*$--equivalences to  $L_{T(n)}\Sinfty(\Oinfty X)_+$ will canonically factor through $s_n$. \\
 
Homological consequences are most precise when $E_* = K(n)_*$.  The first property then says that, for all $X$,
$$ s_n(Z)_*: K(n)_*(\PP(X)) \ra K(n)_*(\Oinfty X)$$
is an inclusion of commutative, cocommutative, $K(n)_*$--Hopf algebras.  The second property and work of Hopkins, Ravenel, and Wilson \cite{hrw} combine to say that, if $X$ is an $S$--module with $T(i)_*(X) = 0$ for $1 \leq i \leq n-1$, then there is an isomorphism of $K(n)_*$--Hopf algebras
$$ K(n)_*(\Oinfty X) \simeq K(n)_*(\PP X) \otimes \bigotimes_{j=0}^{n+1} K(n)_*(K(\pi_j(X),j)).$$

Our main theorems also have consequences for $E_n^*$, where $E_n$ is fundamental $p$--complete integral height $n$ complex oriented commutative $S$--algebra appearing in the work of Hopkins and his collaborators, since it is known \cite{hovey1} that $K(n)_*(X) = 0$ if and only if $E_n^*(X)=0$.  Indeed our work here, combined with work by many, beginning with \cite{hkr} and \cite{hunton}, on $E^*_n(D_rX)$, gets us most of the way towards calculations of $E_n^*(\Oinfty X)$ generalizing Bousfield's extensive functorial calculations  \cite{bousfield4} of $E_1^*(\Oinfty X) = K^*(\Oinfty X; \mathbb Z_p)$.

Our theory of localized Andr\'{e}--Quillen towers also yields the following theorem, a significant generalization of the main result of \cite{k1}: if $Z$ is a connected space, and $f: \Sinfty Z \ra X$ is an $E_*$--isomorphism, then $(\Oinfty f)_*: E_*(\Oinfty \Sinfty Z) \ra E_*(\Oinfty X)$ is monic.

We describe our results in detail in the next section.  

Versions of our main theorem, \thmref{Tn theorem}, date from 2000, and have been reported on in various seminar and conference talks since then in both the United States and Europe. \\

\noindent {\bf Acknowledgements} 

Many people deserve thanks for helping me with aspects of this work.  

I thank Greg Arone and Mike Mandell for many tutorials, taught from complementary prospectives, on  Andr\'e--Quillen--Goodwillie towers of various sorts.   What I have learned from them has also been enhanced by the ideas of Randy McCarthy and his students, and by perusing recent versions of \cite{goodwillie3} that Tom Goodwillie kindly supplied me. Some of the conversations with Mike Mandell were during a visit to the University of Chicago during the fall of 2000, and I thank the Chicago Mathematics Department for its support.  I thank Charles Rezk for helping me through a point of confusion, and Steve Wilson for helping me reach a point of clarity, both related to proofs of results in \secref{sn iso subsection}.

Finally I need to thank Pete Bousfield.   He has been a constant guide to my understanding of the strange behavior of periodic localization.  This has been true for twenty years, but email exchanges dating from mid 2000 have particularly helped me keep straight the technical details of this project.  I particularly recommend his recent paper \cite{bousfield6} for discussions of problems closely related to those studied here.

\section{Main results}

In this section, we describe our results. It ends with a discussion of the organization of the remainder of the paper, where proofs and more detail are given.

Throughout we will use the following convention: if $A$ and $B$ are objects in a model category $\C$, a {\em weak} map $A \xra{f} B$ will mean either a pair $A \xla[\sim]{g} C \xra{h} B$, or a pair $A \xra{h} C \xla[\sim]{g} B$.  A weak map in $\C$ induces a well defined morphism in the homotopy category $ho(\C)$, and we say that a diagram of weak maps in $\C$ commutes if the induced diagram in $ho(\C)$ does.

\subsection{Commutative $S$--algebras and the stable splitting of $QZ$, revisited}

Let $Z_+$ denote the union of a space $Z$ with a disjoint basepoint. If $X$ is an $S$--module, $\Sinfty (\Oinfty X)_+$ is naturally a commutative $S$--algebra augmented over the sphere spectrum $S$.  We denote by $\Alg$ the category of such objects. Another example is $\PP X$, the free commutative algebra generated by $X$.  There is a natural weak equivalence:
\begin{equation*}  \PP X \simeq \bigvee_{r = 0}^{\infty} D_rX.
\end{equation*}

Given $A \in \Alg$, let $I(A)$ be the homotopy fiber of the augmentation $A \ra S$. We view $I(A)$ as the augmentation ideal, and we are interested in two associated objects.  The first is $\widehat A \in \Alg$, which arises as the inverse limit of an {\em Andr\'e--Quillen tower} in $\Alg$, and can be viewed as the $I(A)$--adic completion of $A$.  The second is an $S$--module $taq(A)$, a form of Topological Andr\'{e}--Quillen homology, and can be viewed as $I(A)/I(A)^2$.  There is a convergence result: if $A$ is 0--connected then the canonical map $A \ra \widehat{A}$ is an equivalence.  

Applied to the examples above, if $A$ is either $\PP X$ or $\Sinfty (\Oinfty X)_+$, with $X$ $-1$--connected in the latter case, then $taq(A) \simeq X$.  The natural map $I(A) \ra taq(A)$ identifies in the first case with projection onto the first factor, and in the second with $\epsilon(X): \Sinfty \Oinfty X \ra X$, the counit of the adjunction.  There is also a natural weak equivalence:
\begin{equation*}  \widehat{\PP} X \simeq \prod_{r = 0}^{\infty} D_rX.
\end{equation*}

We observe: 

\begin{prop} \label{tower prop1} Let $f:A \ra B$ be a map in $\Alg$.  If $taq(f): taq(A) \ra taq(B)$ is an equivalence, so is $\widehat{f}: \widehat{A} \ra \widehat{B}$.  Thus, in this case, there is a factorization by weak algebra maps
\begin{equation*}
\xymatrix{
A \ar[rr]^{canonical} \ar[dr]_{f}  &  &  \widehat{A}, \\
   & B \ar[ur] &    \\
}
\end{equation*}
where the unlabelled weak map is $B \ra \widehat{B} \xla[\sim]{\widehat{f}} \widehat{A}$.
\end{prop}

Observations like this are the basis for various of our splitting theorems.  We first illustrate this idea by giving a new formulation and proof of a very highly structured version of the classical splitting of $\Sinfty Q Z$, where, as usual,  $QZ$ denotes $\Oinfty \Sinfty Z$. (See \cite[Appendix B]{k3} for a discussion of some of the different proofs this theorem of D.S.Kahn.)

Let $\eta(Z): Z \ra QZ$ be the unit of the adjunction.   In a straightforward way, this then induces a weak map in $\Alg$ 
$$s(Z): \PP(\Sinfty Z) \ra \Sinfty(QZ)_+.$$

\begin{thm}  \label{classic theorem} For all spaces $Z$, the map $s(Z)$ induces an isomorphism on completions, and thus there is a natural factorization of weak algebra maps
\begin{equation*}
\xymatrix{
\PP(\Sinfty Z) \ar[rr]^{canonical} \ar[dr]_{s(Z)}  &  &  \widehat{\PP}(\Sinfty Z). \\
   & \Sinfty (QZ)_+ \ar[ur]_{t(Z)} &    \\
}
\end{equation*}
If $Z$ is 0--connected then all of these maps are weak equivalences.
\end{thm}

In more down--to--earth terms, our theorem gives us a factorization by weak maps
\begin{equation*}
\xymatrix{
\bigvee_{r = 0}^{\infty} \Sinfty D_rZ \ar[rr]^{canonical} \ar[dr]_{s(Z)}  &  &  \prod_{r=0}^{\infty} \Sinfty D_rZ \\
   & \Sinfty (QZ)_+ \ar[ur]_{t(Z)} &    \\
}
\end{equation*}
in which both the infinite wedge and product are equivalent to $E_{\infty}$--ring spectra, such that all maps in the diagram are $E_{\infty}$.

Our general theory leads to a short proof of the theorem as follows. The diagram
\begin{equation} \label{space lifting}
\xymatrix{   & \Sinfty \Oinfty \Sinfty Z \ar[dr]^{\epsilon(\Sinfty Z)} &    \\
\Sinfty Z \ar@{=}[rr] \ar[ur]^{\Sinfty \eta(Z)}  &  &  \Sinfty Z \\
}
\end{equation}
commutes for all $Z$.  By construction, this shows that $taq(s(Z))$ can be identified with the identity map on $\Sinfty Z$ and thus is an equivalence.  By the proposition, so also is $\widehat{s(Z)}$, and the theorem follows. \\

In \appref{splitting appendix}, we check that our stable splitting agrees with others in the literature. 

In the case when $Z$ is not connected, our theorem improves upon various weaker versions in the literature, and our proof shows that most of the technical issues confronted in these papers need no longer be part of the story.  See Remarks \ref{classic thm remarks}.

\subsection{Bousfield localization and Andr\'{e}--Quillen towers}

Now we mix Bousfield localization with the general theory.  If $E$ is any $S$--module, $L_ES$ will be a commutative $S$--algebra, and we define $L_E(\Alg)$ to be the category of commutative $L_ES$--algebras, which are also augmented over $L_ES$, and are $E$--local.  Up to weak equivalence, objects have the form $L_EA$, with $A\in \Alg$, but {\em not} all morphisms are homotopic to one of the form $L_Ef$, with $f \in \Alg$.

Analogous to the nonlocalized theory, given $A \in L_E(\Alg)$, one gets a completion $\widehat{L}_E A \in L_E(\Alg)$, and an $L_ES$--module $taq^E(A)$.  For $A \in \Alg$,  $L_E(taq(A))\simeq taq^E(L_E A) $, but it is {\em not} in general true that the natural map $L_E \widehat{A} \ra \widehat{L}_E A$
is an equivalence. This is illustrated by the example 
$$ \widehat{L}_E \PP X \simeq \prod_{r=0}^{\infty} L_E D_r X,$$
which is often different than 
$$ L_E\widehat{\PP}X \simeq L_E \left(\prod_{r=0}^{\infty}D_r X\right).$$

Analogous to \propref{tower prop1}, we observe:

\begin{prop} \label{LE tower prop1} Let $g:L_E A \ra L_E B$ be a map in $L_E(\Alg)$.  If $taq^E(g): L_Etaq(A) \ra L_Etaq(B)$ is an equivalence, so is $\widehat{g}: \widehat{L}_E A \ra \widehat{L}_E B$.  Thus, in this case, one gets a factorization by weak algebra maps
\begin{equation*}
\xymatrix{
L_EA \ar[rr]^{canonical} \ar[dr]_{g}  &  &  \widehat{L}_E A. \\
   & L_E B \ar[ur] &    \\
}
\end{equation*}
\end{prop}

In the early 1980's, the author proved that if $f: \Sinfty Z \ra \Sinfty W$ is an $E_*$--isomorphism, with $Z$ and $W$ connected, then $\Oinfty f$ is an $E_*$--monomorphism \cite{k1}.  As a first application of our general theory of localized towers, we  deduce the following stronger version. \\

\begin{thm} \label{susp thm}  Let $Z$ be a connected space.  If a map of spectra $f: \Sinfty Z \ra X$ is an $E_*$--isomorphism, then 
$$ (\Oinfty f)_*: E_*(Q Z) \ra E_*(\Oinfty X)$$
is a monomorphism. \\
\end{thm}

We leave to the reader the proof that the hypotheses can be weakened slightly: the domain of $f$ need just be `spacelike', i.e. a wedge summand of a suspension spectrum.  See also \appref{splitting appendix} for a version of the theorem for non--connected $Z$. 

Examples illustrating this theorem were given in \cite{k1}.  Besides these, see also its use in \appref{Bousfield appendix}.  \\

\subsection{The main theorem}
Fixing a prime $p$ and $n \geq 1$, we now apply the proposition of the last subsection to the case when $E=T(n)$.  Our theorem has a statement and proof analogous to \thmref{classic theorem}.

In place of (\ref{space lifting}), we use the following much deeper theorem: 
there is a natural factorization by weak $S$--module maps
\begin{equation} \label{Tn lifting}
\xymatrix{   & L_{T(n)}\Sinfty \Oinfty X \ar[dr]^{L_{T(n)}\epsilon(X)} &    \\
L_{T(n)}X \ar@{=}[rr] \ar[ur]^{\eta_n(X)}  &  &  L_{T(n)}X. \\
}
\end{equation}
Such a factorization was constructed in the mid 1980's by Bousfield \cite{bousfield3}, when $n=1$, and by the author, for all $n \geq 1$ \cite{k2}.  A recent paper by Bousfield \cite{bousfield6} revisits these constructions.  These papers heavily use the classification of stable $v_n$--self maps of finite complexes, and thus are using (at least for $n > 1$) the work of Devanitz, Hopkins, and Smith \cite{dhs, hs} on Ravenel's Nilpotence Conjectures.

The natural transformation $\eta_n(X)$ then induces a weak map 
$$s_n(X): L_{T(n)}\PP(X) \ra L_{T(n)}\Sinfty (\Oinfty X)_+$$
in $L_{T(n)}(\Alg)$. \propref{LE tower prop1} combines with  (\ref{Tn lifting}) to prove the main theorem of the paper:

\begin{thm}  \label{Tn theorem} For all spectra $X$, the map $s_n(X)$ induces an isomorphism on $L_{T(n)}$--completions, and thus there is a natural factorization of weak algebra maps
\begin{equation*}
\xymatrix{
L_{T(n)}\PP(X) \ar[rr]^{canonical} \ar[dr]_{s_n(X)}  &  &  \widehat{L}_{T(n)}\PP(X). \\
   & L_{T(n)}\Sinfty (\Oinfty X)_+ \ar[ur]_{t_n(Z)} &    \\
}
\end{equation*}
\end{thm}

In more down--to--earth terms, our theorem gives us a factorization by weak maps
\begin{equation*}
\xymatrix{
L_{T(n)} \left( \bigvee_{r = 0}^{\infty} D_rX \right) \ar[rr]^{canonical} \ar[dr]_{s_n(X)}  &  &  \prod_{r=0}^{\infty} L_{T(n)} D_rX \\
   & L_{T(n)} \Sinfty (\Oinfty X)_+ \ar[ur]_{t_n(X)} &    \\
}
\end{equation*}
in which both the infinite wedge and product are equivalent to commutative $L_{T(n)}S$--algebras, such that all maps in the diagram are $L_{T(n)}S$--algebra maps.  \\

For applications to computing $K(n)_*$ and $E_n^*$, it is useful to let $$s_n^K(X) = L_{K(n)}s_n(X): L_{K(n)}\PP(X) \ra L_{K(n)}\Sinfty (\Oinfty X)_+.$$  

The functor on $S$--modules sending $X$ to $\PP(X)$ preserves $E_*$--isomorphisms for any generalized homology theory $E_*$.  Our natural transformations $s_n(X)$ and $s_n^K(X)$ yield the {\em best} possible `invariant' approximations to the functors $L_{T(n)}\Sinfty (\Oinfty X)_+$ and $L_{K(n)}\Sinfty (\Oinfty X)_+$ in the following sense: \\

\begin{prop} \label{sn is universal prop}  Let $F: \Sp \ra \Sp$ be any functor preserving \\
$T(n)_*$--isomorphisms.  Then any natural transformation $T$ of the form
$$ T(X): F(X) \ra L_{T(n)}\Sinfty (\Oinfty X)_+$$
factors uniquely through $s_n$.  Similarly, $s_n^K$ is the terminal natural transformation from a functor preserving $K(n)_*$--isomorphisms. \\
\end{prop}

This proposition will be an easy consequence of results described in \secref{sn iso subsection}. \\

\subsection{First homological corollaries} 

It is easily verified that the canonical map from a wedge to the product of a family of spectra induces a monomorphism on any homology theory.  Thus \thmref{Tn theorem} has the following theorem as an immediate corollary. \\

\begin{thm} If $E_*$ is any homology theory such that $T(n)_*$--acyclics are also $E_*$--acyclics, then $s_n$ induces a natural monomorphism
$$ s_n(X)_*: E_*(\PP X) \ra E_*(\Oinfty X).$$
\end{thm} 

The commutative H--space structure on $\Oinfty X$, together with the diagonal, induces a $K(n)_*$--Hopf algebra structure on $K(n)_*(\Oinfty X)$.\footnote{In the twisted sense described in \cite[Appendix]{bousfield4}, if $p=2$.}  Meanwhile, the usual product maps, 
$$ D_iX \sm D_jX \ra D_{i+j}X,$$
and transfer maps,
$$ D_{i+j}X \ra D_iX \sm D_jX ,$$
associated to the inclusion of groups $\Sigma_i \times \Sigma_j \subseteq \Sigma_{i+j}$ make $K(n)_*(\PP X)$ into a $K(n)_*$--Hopf algebra.   When specialized to $E_*=K(n)_*$, the last theorem refines as follows.

\begin{thm} \label{Kn thm} $s_n(X)_*: K(n)_*(\PP X) \ra K(n)_*(\Oinfty X)$ 
is a natural inclusion of commutative, cocommutative $K(n)_*$--Hopf algebras. 
\end{thm}

Since the cohomological Bousfield class of $E_n$ is the same as the homological Bousfield class of $K(n)$, \thmref{Tn theorem} also has the following consequence.

\begin{thm} \label{En thm}
For all spectra $X$, and all $n \geq 1$, there is a factorization of commutative $E_n^*$--algebras 
\begin{equation*}
\xymatrix{
\bigoplus_{r = 0}^{\infty} E_n^*(D_rX)  \ar[rr]^{canonical} \ar[dr]_{t_n(X)^*}  &  &  \prod_{r=0}^{\infty} E_n^*(D_rX), \\
   & E_n^*(\Oinfty X) \ar[ur]_{s_n(X)^*} &    \\
}
\end{equation*}
where the algebra structure on the infinite sum and product is induced by the transfer maps.  In particular, $t_n(X)^*$ is a natural inclusion of commutative $E_n^*$--algebras.
\end{thm}

Here the map labelled $t_n(X)^*$ is defined by applying $E_n^*$ to the composite
$$ \Sinfty (\Oinfty X)_+ \xra{t_n(X)} \prod_{q=0}^{\infty} D_qX \ra \prod_{q=0}^{r} D_qX$$
and then letting $r$ go to $\infty$.

We explain the appearance of the transfer maps, e.g. in this last theorem.  This is a consequence of the naturality of $s_n$ and $t_n$, as applied to the diagonal $\Delta: X \ra X \times X$.   Since $\Oinfty$ commutes with products, one sees that the product on $E_n^*(\Oinfty X)$ is induced by applying the functor $\Sinfty (\Oinfty ( \hspace{.1in})_+)$ to $\Delta$.
Meanwhile, it is well known (see \cite[Thm.VII.1.10]{lmms} or \cite[Prop.A.3]{kuhn5}) that applying $D_k(\hspace{.1in})$ to the weak map
$$ X \xra{\Delta} X \times X \xla{\sim} X \vee X$$
yields the product, over $i+j = k$, of the transfer maps
$$ D_kX \ra D_iX \sm D_jX.$$

\subsection{When is $s_n(X)$ an equivalence?}  \label{sn iso subsection}

One might now wonder how often $s_n(X)$ and $s_n^K(X)$ are equivalences.   We have various results which together give a good sense of what is happening.

We define classes of $S$--modules $\Sp_n \subseteq \Sp_n^K$ as follows.
\begin{equation*}
\begin{split}
\text{Let } \Sp_n & 
= \{ X \in \Sp \ | \ s_n(X) \text{ is an equivalence} \} \\
  & = \{ X \in \Sp \ | \ s_n(X)_*: T(n)_*(\PP X) \xra{\sim} T(n)_*(\Oinfty X) \}.
\end{split} 
\end{equation*}
\begin{equation*}
\begin{split}
\text{Let } \Sp_n^K & 
= \{ X \in \Sp \ | \ s_n^K(X) \text{ is an equivalence} \} \\
  & = \{ X \in \Sp \ | \ s_n(X)_*: K(n)_*(\PP X) \xra{\sim} K(n)_*(\Oinfty X)  \}.
\end{split} 
\end{equation*}

We recall that Eilenberg--MacLane spectra are acyclic in $T(n)_*$.  From this the following first observations are easily deduced: if $\Sp^{\prime}$ is either $\Sp_n$ or $\Sp_n^K$, then $X \in \Sp^{\prime}$ if and only if $X \langle -1 \rangle \in \Sp^{\prime}$, and furthermore, a necessary condition is that $\pi_0(X) = 0$.  Here $X\langle d \rangle$ denotes the $d$--connected cover of an $S$--module $X$.  

Our results about $\Sp_n$ and $\Sp_n^K$ are most pleasantly described by first introducing two more classes of $S$--modules.
\begin{equation*}
\begin{split}
\text{Let } \bar{\Sp}_n & 
= \{ X \in \Sp \ | X\langle d \rangle \in \Sp_n \text{ for large } d \}.
\end{split} 
\end{equation*}
\begin{equation*}
\begin{split}
\text{Let } \bar{\Sp}_n^K & 
= \{ X \in \Sp \ | \ X\langle d \rangle \in \Sp_n^K \text{ for large } d \}.
\end{split} 
\end{equation*}

We recall a concept from \cite{hrw}: say $X$ is {\em strongly $E_*$--acyclic} if the spaces $\Oinfty \Sigma^c (X\langle -1 \rangle)$ are $E_*$--acyclic for all large $c$.\footnote{Note that $\Oinfty \Sigma^c (X\langle -1 \rangle)$ is the $c^{th}$ space in the connective cover of the spectrum $X$.}   For example, Eilenberg--MacLane spectra are strongly $T(n)_*$--acyclic.

Following notation in various papers, e.g. \cite{bousfield4}, let $L_{n-1}^f$ denote localization with respect to $T(0) \vee \dots \vee T(n-1)$.   We recall that this is smashing: $L_{n-1}^f X \simeq L_{n-1}^f S \sm X$.

Armed with this terminology, we have the following theorem and proposition. \\

\begin{thm} \label{sn bar iso thm}  \hspace{0in} 
 
\noindent (1) $ \bar{\Sp}_n = \{ X \in \Sp \ | \ L_{n-1}^f X \text{ is strongly $T(n)_*$--acyclic} \}$, and  
$\Sp_n \subset \bar{\Sp}_n.$ 

\noindent (2) $ \bar{\Sp}^K_n = \{ X \in \Sp \ | \ L_{n-1}^f X \text{ is strongly $K(n)_*$--acyclic} \}$, and $\Sp_n^K \subset \bar{\Sp}^K_n.$  \\
\end{thm}

\begin{prop}  \label{Lnf strongly acyclic prop} There are implications $(1) \Rightarrow (2) \Rightarrow (3) \Rightarrow (4) \Rightarrow (5)$. 

\noindent (1) $T(i)_*(X) = 0$ for $1 \leq i \leq n-1$. 

\noindent (2) $L_{n-1}^f X$ is strongly $T(n)_*$--acyclic. 

\noindent (3) $L_{n-1}^f X$ is strongly $K(n)_*$--acyclic. 

\noindent (4) (with $n \geq 2$) $X$ is strongly $K(n-1)_*$--acyclic. 

\noindent (5) $K(i)_*(X) = 0$ for $1 \leq i \leq n-1$. \\
\end{prop}

Let $c(n)$ denote the smallest integer $c$ such that $\tilde{T}(n)_*(K(\mathbb Z/p,c)) = 0$.  Then $c(n) \geq n+1$, with equality certainly holding if the Telescope Conjecture is true, and perhaps even if not. \\

\begin{thm} \label{sn(X) iso thm}  \hspace{0in}

\noindent (1) Suppose $X \in \bar{\Sp}_n$.  Then $X \in \Sp_n$ if $\pi_0(X) = 0$, $\pi_{j}(X)$ is uniquely $p$--divisible for $0 \leq j \leq c(n)$, and also $\pi_{c(n)+1}(X)/(\text{torsion})$ is uniquely $p$--divisible. 

\noindent (2) Suppose $X \in \bar{\Sp}_n^K$.  Then $X \in \Sp_n^K$ if and only if $\pi_0(X) = 0$, $\pi_{j}(X)$ is uniquely $p$--divisible for $1 \leq j \leq n$, and also $\pi_{n+1}(X)/(\text{torsion})$ is uniquely $p$--divisible. 
 \\
\end{thm}

\begin{ex}
The conditions in \propref{Lnf strongly acyclic prop} trivially hold for all $X$ if $n=1$. We conclude that that $\bar{\Sp}_1 = \bar{\Sp}_1^K = \Sp$, and $\Sp_1 = \Sp_1^K$ is determined by the second statement  of \thmref{sn(X) iso thm}.  In particular, there is a natural equivalence of commutative augmented $L_{K(1)}S$--algebras
$$ L_{K(1)} \PP(X) \simeq L_{K(1)}\Sinfty(\Oinfty X)_+$$
for all 1--connected $X$, with torsion $\pi_2$.  This fits perfectly with the many results on $K^*(\Oinfty X)$ proved by Bousfield beginning with \cite{bousfield1}. \\
\end{ex}

\begin{ex}
 The Telescope Conjecture holds for $n=1$: $T(1)_*$--acyclics are $K(1)_*$--acyclics.  Thus the conditions in \propref{Lnf strongly acyclic prop} are all equivalent when $n=2$. We conclude that  $\bar{\Sp}_2 = \bar{\Sp}_2^K = \{ X \ | \ K(1)_*(X) = 0 \}$, and $\Sp_2^K$ is the set of $S$--modules described by the second statement of \thmref{sn(X) iso thm}.  In particular, there is a natural equivalence of commutative augmented $L_{K(2)}S$--algebras
$$ L_{K(2)} \PP(X) \simeq L_{K(2)}\Sinfty(\Oinfty X)_+$$
for all $K(1)_*$--acyclic, 2--connected $X$, with torsion $\pi_3$. \\
\end{ex}

\begin{ex}  These results tell us precisely which $p$--local finite spectra are in $\Sp_n^K$ and $\bar{\Sp}_n^K$.  All finites are in $\bar{\Sp}_1^K$, and $F$ is in $\Sp_1^K$ if and only if $\pi_0(F) = \pi_1(F) = 0$ and $\pi_2(F)$ is torsion.  If $n \geq 2$, $F$ is in $\bar{\Sp}_n^K$ if and only if $F$ has type at least $n$, and $F$ is in $\Sp_n^K$ if and only if $F$ has type at least $n$ and also $\pi_j(F) = 0$ for all $0 \leq j \leq n$.  \\
\end{ex}

\begin{rems}  If the Telescope Conjecture is true for a pair $(p,n)$, then the second statement of \thmref{sn(X) iso thm} improves the first by one dimension.  This is due to the fact that in proving this second statement, we use computational methods based on special properties of $K(n)_*$.

Even if the Telescope Conjecture fails, it is still conceivable that some of the conditions in \propref{Lnf strongly acyclic prop} are equivalent.  We note that, as observed in \cite[Thm.3.14]{hrw}, if $X$ is a $BP$--module, then condition (5) implies condition (1). 

Our results imply that, in general, the two maps
$$s_n^K(\Sinfty Z), L_{K(n)}s(Z): L_{T(n)} \PP(\Sinfty Z) \ra L_{T(n)} \Sinfty QZ_+$$ are {\em not} homotopic, since, by \thmref{classic theorem}, the second map is an equivalence whenever $Z$ is connected, while the first map needn't be.  By perturbing the `usual' stable splitting of $QZ$, we have thus lost homotopy equivalence but gained naturality with respect to {\em stable} maps between suspension spectra.  See \appref{hopf invariant appendix} for more about this.
\end{rems}

\subsection{More homological corollaries} \label{more homo cor subsection}

Hand in hand with the results of the last subsection, are some more homological corollaries.

In the spirit of \cite[\S 11]{bousfield7}, if Eilenberg--MacLane spectra are strongly $E_*$--acyclic, one can define $E_*^{vir}(Z)$, the {\em virtual} $E_*$--homology of a space $Z$, by the formula
$$ E_*^{vir}(Z) = E_*(Z \langle d \rangle) \text{ for large } d.$$

Methods of \cite{hrw} imply the following illuminating lemma.

\begin{lem}  \label{strongly acyclic lemma} An $S$--module $X$ is strongly $E_*$--acyclic if and only if $$\Tilde{E}_*^{vir}(\Oinfty X) = 0.$$
\end{lem}

The fact that $\PP (X\langle d \rangle) \ra \PP X$ is a $K(n)_*$--equivalence implies that there is a canonical lifting
\begin{equation*}
\xymatrix{  &
K(n)_*^{vir}(\Oinfty X) \ar[d]  \\
K(n)_*(\PP X) \ar[r]^{s_n(X)_*} \ar[ur]&
K(n)_*(\Oinfty X).
}
\end{equation*}

The next theorem is closely related to \thmref{sn bar iso thm}. \\

\begin{thm}  \label{sn* ses thm} For all $X$, there is a natural short exact sequence of commutative, cocommutative $K(n)_*$--Hopf algebras
$$ K(n)_*(\PP X) \xra{s_n(X)_*} K(n)_*^{vir}(\Oinfty X) \ra K(n)_*^{vir}(\Oinfty L_{n-1}^f X).$$
\end{thm}

Note that the first term here is a functor of $L_{K(n)}X$; in constrast, the last term is `invisible' to $L_{K(n)}X$, as $K(n)_*(L_{n-1}^f X) = 0$.  

Since $ K(n)_*(\PP X)=0$ exactly when $K(n)_*(X) = 0$, we have the next corollary, which strengthens \cite[Cor.3.13]{hrw}. \\

\begin{cor} \label{Kn strongly acyclic cor}  $X$ is strongly $K(n)_*$--acyclic if and only if $X$ is $K(n)_*$--acyclic and also $L_{n-1}^f X$ is strongly $K(n)_*$--acyclic.  In particular, if $X$ is $K(n)_*$--acyclic, and also $T(i)_*$--acyclic for $1 \leq i\leq n-1$, then $X$ is strongly $K(n)_*$--acyclic.  \\
\end{cor}

Our next result is a homological variant of  \propref{sn is universal prop}. \\

\begin{prop} \label{sn* is universal prop}  Let $F: \Sp \ra K(n)_*$--modules be any functor preserving $K(n)_*$--isomorphisms.  Then any natural transformation $T$ of the form
$$ T(X): F(X) \ra K(n)_*(\Oinfty X)$$
factors uniquely through $s_{n*}$.  \\
\end{prop}

The $T(n)_*$ variant of this proposition also holds.  Similarly, there is also a $E_n^*$ variant that says that $s_n^*$ is the initial functor from $E_n^*(\Oinfty X)$ to a functor preserving $E_n^*$--isomorphisms. \\

The next theorem has \thmref{sn(X) iso thm}(2) as a consequence. \\

\begin{thm} \label{Kn thm 2}  Let $X \in \bar{\Sp}_n^K$, i.e. $L_{n-1}^f X$ is strongly $K(n)_*$--acyclic.  Then each of the maps
$$ K(n)_*(\Oinfty X\langle j \rangle) \ra K(n)_*(\Oinfty X)$$
is an inclusion of a normal sub--$K(n)_*$--Hopf algebra.  This induces a decreasing filtration of finite length on $K(n)_*(\Oinfty X)$, and there  
is an isomorphism of filtered $K(n)_*$--Hopf algebras,
$$ K(n)_*(\Oinfty X) \simeq K(n)_*(\PP X) \otimes \bigotimes_{j=0}^{n+1} K(n)_*(K(\pi_j(X),j)),$$
that is natural on the level of associated graded objects. \\
\end{thm}

\begin{ex} When $n=1$, the theorem says that for all $X$, there is an isomorphism of $K(1)_*$--Hopf algebras
$$ K(1)_*(\Oinfty X) \simeq K(1)_*(\PP X) \otimes \bigotimes_{j=0}^{2} K(1)_*(K(\pi_j(X),j)).$$
\end{ex}

\begin{ex} When $n=2$, the theorem says that, if $K(1)_*(X) = 0$, then there is an isomorphism of $K(2)_*$--Hopf algebras
$$ K(2)_*(\Oinfty X) \simeq K(2)_*(\PP X) \otimes \bigotimes_{j=0}^{3} K(2)_*(K(\pi_j(X),j)).$$
\end{ex}

\begin{ex}   Suppose $L_{n-1}^f X$ is strongly $K(n)_*$--acyclic, and also $X$ is $n$--connected with $\pi_{n+1}(X)$ torsion. The theorem implies that
\begin{equation} \label{conn cover thm}
K(n)_*(\Oinfty X\langle n+1 \rangle) \ra K(n)_*(\Oinfty X)
\end{equation}
is an isomorphism. 

When $n=1$, the first hypothesis is always satisfied, and we recover \cite[Thm. 2.4]{bousfield1}.  This is the key technical theorem of Bousfield's paper.  In recent email to the author, Bousfield has observed that (\ref{conn cover thm}) allows for an addendum to \cite[Thm. 8.1]{bousfield6}, analogous to the use of \cite[Thm. 2.4]{bousfield1} in the proof of \cite[Thm. 3.2]{bousfield1}.

Some hypotheses are necessary here.  In \cite[\S 8.7]{bousfield6}, Bousfield notes that if $X$ is the suspension spectrum of the Moore space $M(\Z/p,3)$, then 
$$ K(2)_*(\Oinfty X\langle 3 \rangle) \ra K(2)_*(\Oinfty X)$$
has nonzero kernel.  \\
\end{ex}

\begin{ex}  Let $k(n)_c = \Oinfty \Sigma^c k(n)$, where $k(n)$ is the connective cover of $K(n)$.  $k(n) \in \bar{\Sp}_n^K$, as it is a $BP$--module that is $K(i)_*$--acyclic for $1 \leq i \leq n-1$.  Thus the theorem applies to say that the cofibration sequence of spectra
$$ \Sigma^{2p^n -2} k(n) \xra{v} k(n) \ra H\Z/p$$
induces a short exact sequence of commutative, cocommutative $K(n)_*$--Hopf algebras
$$ K(n)_*(k(n)_{2p^n -2 +c}) \ra K(n)_*(k(n)_c) \ra K(n)_*(K(\Z/p,c))$$
for all $c > 0$. 

This is \cite[Thm.1.1]{bkw}. \\
\end{ex}

These last two examples also illustrate our next result.  \\

\begin{thm} \label{ses thm}  Suppose $f: X \ra Y$ is map of 0--connected $S$--modules with cofiber $C$, such that $\PP(f)_*: K(n)_*(\PP X) \ra K(n)_*(\PP Y)$ is monic. (For example, $f$ might be a $K(n)_*$--isomorphism.)  If $X \in \Sp_n^K$ then there is a short exact sequence of commutative, cocommutative  $K(n)_*$--Hopf algebras
$$ K(n)_*(\Oinfty X) \ra K(n)_*(\Oinfty Y) \ra K(n)_*(\Oinfty C).$$
\end{thm}

\begin{rems}  It seems appropriate to comment on other results in the literature that concern homological calculations of the sort we look at here.

In unpublished work, but in the spirit of \cite{strickland}, N. Strickland has observed that $E_n^*(\PP X)$ is a functor of $E_n^*(X)$, when restricted to $X$ such that $E_n^*(X)$ is appropriately `pro--free' as an $E_n^*$--module, and also $E_n^*(X)$ is concentrated in even degrees.  It strikes the author as likely that the second of these hypotheses is unnecessary, since the work of either J. McClure \cite[Chapter IX]{bmms} or Bousfield \cite{bousfield4} shows this to be the case when $n=1$.  

Related to this, it appears that when $K(n)_*(X)$ and $K(n)_*(Y)$ are concentrated in even degrees, if $f: X \ra Y$ is a $K(n)_*$--monomorphism, then so is $\PP(f)$.

In various papers culminating in \cite{kashiwabara}, Kashiwabara computes $E^*(\Oinfty X)$ if $E = BP, E_n, \text{ or } K(n)$, under suitable side hypotheses on $X$, e.g. if $X$ is $-1$--connected with cells only in even dimensions.  His methods are very different than ours; in particular, he heavily uses $BP$--Adams resolutions of his spectra.  Most of his results appear to have only limited naturality, and many of his results are specialized to the case when $X$ is a suspension spectrum.

There is most obvious overlap between our results and those of Bousfield, particularly in \cite{bousfield6}.  It seems that any proof of our main theorem, \thmref{Tn theorem}, will depend crucially on the existence of a Goodwillie tower for $\Sinfty \Oinfty X$.  But once this theorem has been established, many of our other results allow for alternate proofs using his work.  See  \appref{Bousfield appendix}. \\
\end{rems}

\subsection{The main theorem for rational homology}

For completeness, we note that a version of our main theorem holds when $n=0$, i.e., with $L_{H\Q}$ replacing $L_{T(n)}$.

Analogous to (\ref{Tn lifting}), we have the following lemma:  for all $0$--connected spectra $X$, there is a natural factorization by weak $S$--module maps
\begin{equation} \label{HQ lifting}
\xymatrix{   & L_{H\Q}\Sinfty \Oinfty X \ar[dr]^{L_{H\Q}\epsilon(X)} &    \\
L_{H\Q}X \ar@{=}[rr] \ar[ur]^{\eta_0(X)}  &  &  L_{H\Q}X. \\
}
\end{equation}
For $X$ just $-1$--connected, such a factorization also exists, but can not be made to be natural.

$\eta_0(X)$ then induces a weak map $s_0(X): L_{H\Q}\PP(X) \ra L_{H\Q}\Sinfty (\Oinfty X)_+$ in $L_{H\Q}(\Alg)$, natural for $0$--connected $X$.  \propref{LE tower prop1} combines with  the lemma to prove:

\begin{thm}  \label{HQ theorem} For all $-1$--connected spectra $X$, the map $s_0(X)$ induces an isomorphism on $L_{H\Q}$--completions, and thus there is a factorization of weak algebra maps
\begin{equation*}
\xymatrix{
L_{H\Q}\PP(X) \ar[rr]^{canonical} \ar[dr]_{s_0(X)}  &  &  \widehat{L}_{H\Q}\PP(X). \\
   & L_{H\Q}\Sinfty (\Oinfty X)_+ \ar[ur]_{t_0(Z)} &    \\
}
\end{equation*}
This is natural for $0$--connected $X$, and, in this case, all three maps are equivalences. \\
\end{thm}

As a corollary, one recovers the known theorem: for $0$--connected spectra $X$, the Hopf algebra $ \pi_*^S(\Oinfty X) \otimes \Q$ is naturally isomorphic to the graded symmetric algebra primitively generated by $\pi_*(X)\otimes \Q$.

\subsection{Organization of the paper}

In \secref{tower section}, we develop the general theory of Andr\'e--Quillen towers associated to commutative augmented $S$--algebras, leading to proofs of \propref{tower prop1} and \propref{LE tower prop1}.  Included is a subsection summarizing the basic properties of the functor $\Sinfty(\Oinfty X)_+$.

The splitting theorems concerning $\Sinfty (\Oinfty X)_+$, Theorems \ref{classic theorem}, \ref{Tn theorem}, and \ref{HQ theorem}, are then easily proved in \secref{proofs section}, which includes discussion of (\ref{Tn lifting}) and (\ref{HQ lifting}), and \thmref{susp thm}.

In \secref{snX section}, we explore the extent to which $s_n(X)$ is an equivalence, proving the results in \secref{sn iso subsection} and \secref{more homo cor subsection}.  Some of the proofs are a bit long and delicate: we hope we have made them comprehensible.

In \appref{splitting appendix}, we check that our stable splitting of $QZ$ agrees with others in the literature.  In \appref{Bousfield appendix}, we compare our constructions and theorems to those of \cite{bousfield6}.  In \appref{hopf invariant appendix}, we compare $s$ to $s_n$, and make some remarks about James--Hopf invariants.

\section{The Andr\'e--Quillen tower of commutative algebras} \label{tower section}

\subsection{Categories of commutative $S$--algebras}

We work always within the topological model category of $S$--modules as in \cite{ekmm}.  This is a symmetric monoidal category with unit the sphere spectrum $S$, and we recall that associative, commutative, unital $S$--algebras are modern day versions of $E_{\infty}$--algebras.  Given such an algebra $R$, we let $R-\Alg$ denote the category of associative, commutative, unital, augmented $R$--algebras.  When $R=S$ we simplify the notation to $\Alg$.

Given $A \in R-\Alg$, let $I^{\prime}(A)$ denote the fiber of the augmentation $A \ra R$.  As discussed in \cite{basterra}, the functor $I^{\prime}$ takes values in the category $R-\Algg$ of associative, commutative, nonunital $R$--algebras, and is the right adjoint of a Quillen equivalence between the model categories $R-\Alg$ and $R-\Algg$.\footnote{The left adjoint sends a nonunital $R$--algebra $I$ to the augmented algebra $I^{\prime} \vee R$.}  We let $I(A)$ denote a cofibrant replacement in $R-\Algg$ of $I^{\prime}(A^{\prime})$, where $A^{\prime}$ is a fibrant replacement of $A \in \Alg$. 

The categories $R-\Alg$ and $R-\Algg$ are tensored and cotensored over based topological spaces. (See \cite{ekmm} or \cite{k4}.)  In particular, given $I \in R-\Algg$, one can form the iterated suspension $S^n \otimes I$ and the iterated looping $\Omega^n I$.

We note that both $\Omega^n$ (as a right adjoint) and homotopy colimits over directed systems (see, e.g. \cite[Lemma VII.3.10]{ekmm}) commute with the forgetful functor from $R-\Algg$ to $R$--modules.  

The processes of taking coproducts and suspending (tensoring with $S^1$) certainly {\em don't} commute with this forgetful functor.  Indeed, the coproduct in $R-\Alg$ is the smash product, and thus in $R-\Algg$ one has 
$$ I \amalg J = I\sm J \vee I \vee J.$$
Regarding suspension, one has the following lemma, which serves as the basis for many convergence results. \\

\begin{lem} \label{convergence lemma} If a cofibrant $I \in \Algg$ is $n$--connected, then 
the natural map 
$$ I \ra \Omega(S^1 \otimes I)$$
is $2n+1$--connected.  \\
\end{lem}

We feature two families of examples. \\

\begin{ex}  \label{free example} Given an $S$--module $X$, we let $\PP X$ denote the free commutative $S$--algebra generated by $X$ \cite[p.40]{ekmm}.  If $X$ is cofibrant then there is a natural weak equivalence \cite[p.64]{ekmm}:
\begin{equation*}  \PP X \simeq \bigvee_{r = 0}^{\infty} D_rX.
\end{equation*}
$\PP X$ is naturally augmented, and we let $i: X \ra I(\PP X)$ denote the natural weak map.  Using the freeness of $\PP X$, given any $A \in \Alg$, a weak map of $S$--modules $f: X \ra I(A)$ induces a weak map in $\Alg$, $\tilde{f}: \PP X \ra A$ such that the diagram of weak maps
\begin{equation*}
\xymatrix{
X \ar[d]^{i} \ar[dr]^{f} &  \\
I(\PP X) \ar[r]^{I(\tilde f)} &  I(A)
}
\end{equation*}
commutes.

Given a commutative $S$--algebra $R$, and an $R$--module $Y$, there is an analogous free object $\PP_R Y \in R-\Alg$, satisfying the evident `change of rings' formula
$$ \PP_R (R \sm X) = R \sm \PP X.$$
\end{ex}

\begin{ex}  \label{SLX example} Given an $S$--module $X$, the $E_{\infty}$--structure on the infinite loopspace $\Oinfty X$ implies that $\Sinfty (\Oinfty X)_+$ takes values in $\Alg$.  See \cite[Ex.IV.1.10]{may} and \cite[\S II.4]{ekmm}.  The composite
$$ I(\Sinfty (\Oinfty X)_+) \ra \Sinfty (\Oinfty X)_+ \ra \Sinfty (\Oinfty X)$$
is a weak equivalence for all $X$.  The natural map $X\langle -1 \rangle \ra X$ of $S$--modules induces an equivalence in $\Alg$:
$$ \Sinfty(\Oinfty X\langle -1 \rangle)_+ \xra{\sim} \Sinfty(\Oinfty X)_+.$$
\end{ex}

\subsection{Topological Andr\'e--Quillen homology}

One version of the Topological Andr\'e Quillen Homology of $A \in R-\Alg$ is as the set of homotopy groups of the following construction.

\begin{defn}  Given $A \in R-\Alg$, let $taq_R(A) \in R-\Algg$ be defined by
$$ taq_R(A) = \hocolim_{n \ra \infty} \Omega^n (S^n \otimes I(A)).$$
\end{defn}

\begin{rem}  An alternative construction, more reminscent of the work of Andr\'e and Quillen is to let $taq_R(A) = ZQ(A)$, where
$Q: \Alg \ra R-\text{modules}$ is defined by  
$$Q(A) = I(A)/I(A)^2,$$
and $Z:R-\text{modules} \ra \Algg$ is defined by letting $Z(X)$ be a fibrant replacement of $X$ given trivial multiplication.  This is the construction explored by M.Basterra in \cite{basterra}, and she and M.Mandell have an unpublished proof that these two constructions are equivalent.\footnote{The proof of the main theorem of \cite{basterramccarthy} indicates some of the ideas, as does S.~Schwede's earlier paper \cite{schwede}.}  For our purposes, particularly as in \exref{SLX example: 2} below, the construction we use is most convenient.  \\
\end{rem}

We denote $taq_S(A)$ by $taq(A)$.  The following `change of rings' formula follows easily from the definition. 

\begin{lem} \label{taq change of rings} For all $A \in \Alg$, there is a natural weak equivalence in $R-\Algg$
$$R \sm taq(A)\simeq taq_R(R \sm A).$$
\end{lem}

Another basic property that we will use is the following.

\begin{lem} \label{taq cofib lem} $taq_R$ takes homotopy cofibration sequences in $R-\Alg$ to homotopy cofibration sequences in $R$--modules.
\end{lem}

\begin{proof} As will be elaborated on in \secref{proof of the tower properties}, the functor $taq_R$, defined via stablization as above, will be 1-excisive in the sense of Goodwillie: it will take homotopy pushout squares to homotopy pullback squares.  But in $R$--modules, homotopy pullbacks are homotopy pushouts.
\end{proof}

We now calculate $taq(A)$ for our two key examples.

\begin{ex} \label{free example: 2} Corresponding to \exref{free example}, we claim that there is a natural equivalence 
$$taq(\PP X) \simeq X,$$
where $X$ has trivial multiplication, such the natural map $I(\PP X) \ra taq(\PP X)$ corresponds to the projection
$$ \bigvee_{r = 1}^{\infty} D_rX \ra D_1 X = X.$$
To see this, we make two observations.  

Firstly, for all based spaces $K$ and $S$--modules $X$, there is a natural isomorphism
$$ K \otimes \PP(X) = \PP(K \sm X).$$

Secondly, for all $S$--modules $X$ there are natural equivalences
\begin{equation*}
\hocolim_{n\ra \infty} \Omega^n D_r(\Sigma^n X) \simeq \begin{cases} X & \text{if } r = 0 \\ * & \text{if } r > 0.
\end{cases}
\end{equation*}
(This is clear by a connectivity argument if $X$ is connective, and then note that an arbitrary $S$--module is equivalent to a hocolimit of connective $S$--modules.)

Combining these observations, we compute:
\begin{equation*}
\begin{split}
taq(\PP X) & = \hocolim_{n \ra \infty} \Omega^n (S^n \otimes I(\PP X)) \\
  & = \hocolim_{n \ra \infty} \Omega^n I(\PP \Sigma^n X) \\
  & \simeq \bigvee_{r=1}^{\infty} \hocolim_{n \ra \infty} \Omega^n D_r(\Sigma^n X) \\
  & \simeq X.
\end{split}
\end{equation*}

\end{ex}

\begin{ex} \label{SLX example: 2} If $X$ is a $-1$--connected $S$--module, corresponding to \exref{SLX example}, we claim that there is a natural weak equivalence
$$ taq(\Sinfty (\Oinfty X)_+) \simeq X,$$
such that the natural map $I(\Sinfty (\Oinfty X)_+) \ra taq(\Sinfty (\Oinfty X)_+)$ corresponds to the counit
$$\epsilon: \Sinfty \Oinfty X \ra X.$$  

To see this, we again make a couple of observations.

Firstly, if $Z$ is an $E_{\infty}$--space, let $BZ$ be the associated classifying space \cite{may1}.  As surveyed in \cite{k4}, there are natural weak equivalences
$$ S^1 \otimes I(\Sinfty Z_+) \simeq I(\Sinfty BZ_+)$$
such that the natural map $I(\Sinfty Z_+) \ra \Omega (S^1 \otimes I(\Sinfty Z_+))$ corresponds to $\Sinfty Z \ra \Omega \Sinfty BZ$.

Secondly, given a $-1$--connected $S$--module $X$, let $X_n = \Oinfty \Sigma^n X$.  Then there is a natural equivalence $BX_n \xra{\sim} X_{n+1}$, and a natural weak equivalence
$$ \hocolim_{n \ra \infty} \Omega^n \Sinfty X_n \xra{\sim} X$$
such that the inclusion of $\Sinfty X_0$ into the hocolimit corresponds to $\epsilon$.

Combining these observations, we compute:
\begin{equation*}
\begin{split}
taq(\Sinfty (\Oinfty X)_+) & = \hocolim_{n \ra \infty} \Omega^n (S^n \otimes I(\Sinfty(X_0)_+)) \\
  & \simeq \hocolim_{n \ra \infty} \Omega^n I(\Sinfty (B^n X_0)_+) \\
  & \simeq \hocolim_{n \ra \infty} \Omega^n I(\Sinfty (X_n)_+) \\
  & \simeq \hocolim_{n \ra \infty} \Omega^n \Sinfty X_n \\
  & \simeq X.
\end{split}
\end{equation*}

\end{ex}
\subsection{The Andr\'e--Quillen tower of an augmented $R$--algebra}

It seems that various people have noted the existence of an `Andr\'e--Quillen tower' associated to $A \in R-\Alg$.  Intuitively, this tower is supposed to be the augmentation ideal tower 
$$   \dots \ra A/I^r \ra \dots \ra A/I^2 \ra A/I$$
constructed in a homotopically meaningful way.  The next theorem lists the properties we care about.  

In this theorem, $D_r^R Y$ denotes the $r^{th}$ extended power construction in the category of $R$--modules, i.e. one uses the smash product $\sm_R$.  Note that there is an isomorphism $R \sm D_r(X) = D_r^R(R \sm X)$. \\

\begin{thm} \label{tower theorem}
Given $A \in R-\Alg$, there is a unique natural tower of fibrations in $R-\Alg$ under $A$ 
\begin{equation}\label{algebra tower}
\xymatrix{
&&& \vdots \ar[d] \\
&&& P_{R,2} A\ar[d]^{p_2} \\
&&& P_{R,1} A \ar[d]^{p_1} \\
A \ar[rrr]^{e_0}  \ar[urrr]^{e_1} \ar[uurrr]^{e_2} &&& P_{R,0} A, 
}
\end{equation}
with the following properties. 

\noindent (1) $P_{R,0} A \simeq R$ so that $e_0$ identifies with the augmentation.

\noindent (2) For $r \geq 1$, the fiber of $ p_r: P_{R,r} A \ra P_{R,r-1} A$
is naturally weakly equivalent to $D_r^R (taq_R(A))$.   Furthermore, $I(e_1)$ identifies with the natural map $I(A) \ra taq_R(A)$.

\noindent (3)  Denoting $P_{S,r}(A)$ by $P_r(A)$, there is a change of rings formula: Given $A \in \Alg$ and $R$ a commutative $S$--algebra, there is a natural weak equivalence of towers under $R \sm A$:
$$R \sm P_r(A) \simeq P_{R,r}(R \sm A).$$

\noindent (4) If $I(A)$ is 0--connected, then $e_r$ is $r$--connected. \\
\end{thm}

The uniqueness statement means up to natural weak equivalence.  The weak equivalence in the second property is as $R$--modules.\footnote{Though we won't need nor prove it, if one gives $D_r^R(taq_R(A))$ trivial multiplication, this equivalence is even as objects in $R-\Algg$.}

V.Minasian constructs a tower with these properties in the preprint \cite{minasian}, following along the lines of \cite{basterra}, and using her version of $taq_R(A)$.\footnote{In email with the author, M.Mandell has also sketched this result.} However, the appearance (finally) of a finished version of \cite{goodwillie3} allows the author to feel comfortable with an alternative construction, suggested to him by G.Arone. \\

\begin{defn}  Let the tower $\{P_{R,r}(\hspace{.1in})\}$ denote the Goodwillie tower associated to the identity functor on $R-\Alg$. \\
\end{defn}

In the subsection \secref{proof of the tower properties} we sketch a proof that this tower has the properties stated in the theorem.  Assuming the theorem, we now follow up with some consequences and an example.

\begin{defn}  Let $\displaystyle \widehat{A} = \holim_{r \ra \infty} P_{R,r} A$.
\end{defn}

\begin{cor} If $I(A)$ is 0--connected, then the natural map $A \ra \widehat{A}$ is an equivalence.
\end{cor}

\begin{ex} \label{free example: 3} Corresponding to \exref{free example}, we have 
\begin{equation} \label{tower for PX} P_r (\PP X) \simeq (\bigvee_{q = 0}^{\infty} D_qX)/(\bigvee_{q = r+1}^{\infty} D_qX) 
\end{equation}
so that there is a natural equivalence 
$$ \widehat{\PP} X \simeq \prod_{r = 0}^{\infty} D_rX.$$

One way to see this is to note that both sides of (\ref{tower for PX}), viewed as towers of functors from $S$--modules to $\Alg$, have the correct form to be a Goodwillie tower of the functor sending $X$ to $\PP X$, and thus agree, up to natural weak equivalence.  \\
\end{ex}

The following lemma is well known and much used. \\

\begin{lem} \label{Dr lemma} The functor $D_r^R$ preserves weak equivalences of $R$--modules.
\end{lem}
\begin{proof} Both $\sm_R$ and $(\hspace{.1in})_{h\Sigma_r}$ (homotopy $\Sigma_r$--orbits) preserve weak equivalences.
\end{proof}

\begin{prop} \label{tower prop} Let $f:A \ra B$ be a map in $R-\Alg$.  If $taq_R(f): taq_R(A) \ra taq_R(B)$ is a weak equivalence, so is $P_{R,r}(f): P_{R,r}(A) \ra P_{R,r}(B)$ for all $r$, and thus also
$\widehat{f}: \widehat{A} \ra \widehat{B}$.  Thus, in this case, one gets a factorization by weak $R$--algebra maps
\begin{equation*}
\xymatrix{
A \ar[rr]^{canonical} \ar[dr]_{f}  &  &  \widehat{A}, \\
   & B \ar[ur] &    \\
}
\end{equation*}
where the unlabelled weak map is $B \ra \widehat{B} \xla[\sim]{\widehat{f}} \widehat{A}$.
\end{prop}
\begin{proof} Using the lemma, one proves this by induction up the Andr\'e--Quillen tower.
\end{proof}

This proposition specializes to \propref{tower prop1} when $R=S$. 

\subsection{A summary of the properties of $\PP(X)$ and $\Sinfty(\Oinfty X)_+$}  In this subsection, we use our work thus far to summarize basic properties of $\PP(X)$ and $\Sinfty (\Oinfty X)_+$, viewed as functors from $S$--modules to $\Alg$. \\

\begin{prop} \label{P prop}  The functor $\PP$ satisfies the following properties. 

\noindent (1) $\PP$ takes homotopy colimits in $S$--modules to homotopy colimits in $\Alg$.

\noindent (2) $I(\PP(X)) \ra taq(\PP(X))$ identifies with $\bigvee_{r=1}^{\infty} D_rX \ra D_1X = X$.  Thus the composite
$ X \xra{i} I(\PP(X)) \ra taq(\PP(X))$ is an equivalence.

\noindent (3) $\PP(X) \ra \widehat{\PP}(X)$ identifies with $\bigvee_{r=0}^{\infty} D_r X \ra \prod_{r=0}^{\infty} D_rX$. 
\end{prop}

\begin{proof} The first property follows formally from the fact that $\PP$ is left adjoint to the forgetful functor from $\Alg$ to $S$--modules, since the model category on $\Alg$ is defined so that algebra maps are fibrations or weak equivalences exactly when they are fibrations or weak equivalences when considered as maps of $S$--modules.  The other two properties were established above in \exref{free example: 2} and \exref{free example: 3}. 
\end{proof}

\begin{prop} \label{SLX prop}  The functor $\Sinfty(\Oinfty \text{\hspace{.1in}})_+$ satisfies the following properties. 

\noindent (1) $\Sinfty(\Oinfty \text{\hspace{.1in}})_+$ takes filtered homotopy colimits  in $S$--modules to filtered homotopy colimits in $\Alg$.

\noindent (2) $\Sinfty(\Oinfty \text{\hspace{.1in}})_+$ takes coproducts in $S$--modules to coproducts in $\Alg$. 

\noindent (3) If $X \ra Y \ra Z$ is a cofibration sequence $S$--modules, with $X$ and $Y$ $-1$--connected and $Z$ $0$--connected, then 
$$\Sinfty(\Oinfty X)_+) \ra  \Sinfty(\Oinfty Y)_+) \ra \Sinfty(\Oinfty Z)_+)$$
is a cofibration sequence in $\Alg$.

\noindent (4) $I(\Sinfty(\Oinfty X)_+) \ra taq(\Sinfty(\Oinfty X)_+)$ identifies with $\epsilon: \Sinfty \Oinfty X \ra X\langle -1 \rangle$. 
\end{prop}

\begin{proof} The last property was established above in \exref{SLX example: 2}. 

To see that the first property holds, we first note that filtered homotopy colimits in $\Alg$ are detected by viewing them as being in $S$--modules. (Compare with \cite[\S II.7]{ekmm}.)  But, as a functor to $S$--modules, $\Sinfty(\Oinfty \text{\hspace{.1in}})_+$ certainly commutes with filtered homotopy colimits.

Thanks to the first property, it suffices to prove the second property for finite coproducts.  In $\Alg$, we have equivalences
$$\Sinfty(\Oinfty (X \vee Y))_+ \xra{\sim} \Sinfty(\Oinfty (X \times Y))_+ = \Sinfty(\Oinfty X)_+ \sm \Sinfty(\Oinfty Y)_+, $$ which is the coproduct in $\Alg$ of $\Sinfty(\Oinfty X)_+$ and $\Sinfty(\Oinfty Y)_+$.  

The proof of the third property is more delicate.  A cofiber sequence of $S$--modules
$$ X \xra{f} Y \xra{g} Z$$
will induce a commutative diagram in $\Alg$:
\begin{equation*}
\xymatrix{
\Sinfty(\Oinfty X)_+ \ar@{=}[d] \ar[r] & \Sinfty(\Oinfty Y)_+ \ar@{=}[d] \ar[r] & A \ar[d]^h \\
\Sinfty(\Oinfty X)_+  \ar[r] & \Sinfty(\Oinfty Y)_+  \ar[r] & \Sinfty(\Oinfty Z)_+ \\
}
\end{equation*}
where $A$ is the cofiber in $\Alg$ of $\Sinfty(\Oinfty f)_+$.  We wish to show the map $h$ is an equivalence.

If $X$, $Y$, and $Z$, are all $-1$--connected, then applying $taq$ to this diagram yields the diagram
\begin{equation*}
\xymatrix{
X \ar@{=}[d] \ar[r]^f & Y \ar@{=}[d] \ar[r] & taq(A) \ar[d]^{taq(h)} \\
X  \ar[r]^f & Y \ar[r]^g & Z \\
}
\end{equation*}
where we have used property (4) above.  The bottom horizontal sequence is given as a cofibration sequence; by \lemref{taq cofib lem}, so is the top horizontal sequence.  We conclude that $taq(h)$ is an equivalence.  By \propref{tower prop1}, we conclude that $\widehat{h}$ is an equivalence.  Under our connectivity hypothesis that $Z$ is also $0$--connected (so that $\pi_0(f)$ is onto), one can deduce that both $A$ and $\Sinfty(\Oinfty Z)_+$ have 0--connected augmentation ideals, and thus are equivalent to their completions, by \thmref{tower theorem}(4).  We conclude that $h$ is an equivalence. 
\end{proof}

\begin{rem} If one regards $\Sinfty \Oinfty X_+$ just as a functor taking values in $S$--modules, rather than in $\Alg$, the fact that its Goodwillie tower has $r^{th}$ fiber equivalent to $D_r(X)$ has been known for awhile by Goodwillie and others working with the calculus of functors.  This tower appears explicitly in the literature in \cite{ahearnkuhn}.
\end{rem}

\subsection{The localized tower} If $E$ is an $S$--module, let $L_E$ denote Bousfield localization with respect to $E_*$.  It has long been usefully observed that various constructions in infinite loopspace theory behave well with respect to Bousfield localization.  For example, \cite{k1} heavily used the follow analogue of \lemref{Dr lemma}. 

\begin{lem} \label{Dr Eiso lemma} \cite[Cor.2.3]{k1} The functor $D^R_r$ preserves $E_*$--isomorphisms. 
\end{lem}
\begin{proof} Both $\sm_R$ and $(\hspace{.1in})_{h\Sigma_r}$  preserve $E_*$--isomorphisms.
\end{proof}

This same fact is behind the beautiful and much more recent theorem that if $R$ is a commutative $S$--algebra, so is $L_ER$ \cite[Chap.VIII]{ekmm}. 

\begin{lem} \label{taq Eiso lemma} The functor $taq_R$ preserves $E_*$--isomorphisms. 
\end{lem}
\begin{proof}  \cite[Cor.7.5]{k4} says that $taq_R(A)$ is the colimit of an increasing filtration $F_1taq_R(A) \ra F_2taq_R(A) \ra \dots$ by cofibrations, and identifies the cofibers: there is an equivalence
$$ F_d taq_R(A)/F_{d-1} taq_R(A) \simeq (\Sigma K_d \sm_R I(A)^{\sm_R d})_{h\Sigma_d},$$
where $K_d$ is a certain partition complex appearing in \cite{am}.  The functor on the right of this equivalence certainly preserves $E_*$--isomorphisms, and thus so does $taq_R(A)$. 
\end{proof}

The lemmas combine with induction up the Andr\'e--Quillen tower to prove \\

\begin{cor} \label{Pr Eiso cor} The functors $P_{R,r}$ preserve $E_*$--isomorphisms. \\
\end{cor}

\begin{defn} Let $L_E(\Alg)$ be the full subcategory of $L_ES-\Alg$ consisting of $E$--local objects. \\ 
\end{defn}

In the spirit of these last results, we have the following proposition. \\

\begin{prop}  \label{LE model cat prop}  $L_E: \Alg \ra L_E(\Alg)$ commutes with homotopy pushout squares and filtered homotopy colimits in the following sense: 

\noindent (1)  the natural map
$$ L_E(\hocolim \{B \la A \ra C \}) \ra L_E(\hocolim \{L_EB \la L_EA \ra L_EC \})$$
is an equivalence, and

\noindent (2) the natural map $L_E(\hocolim_{i} A_i) \ra L_E(\hocolim_{i} L_E A_i)$ is an equivalence.
\end{prop}

\begin{proof} As discussed on \cite[p.162]{ekmm}, a model for the pushout of a diagram of $R$--algebras of the form $B \la A \ra C$, where both maps are cofibrations, is given by an appropriate bar construction $\beta^R(B,A,C)$.  This construction preserves $E_*$--isomorphisms in all variables; in particular
$$ L_E(\beta^S(B,A,C)) \ra L_E(\beta^{L_ES}(L_EB,L_EA,L_EC))$$
is an equivalence, establishing (1). The proof of (2) is similar. \\
\end{proof}

\begin{rem} It is easy to see that $B \in L_ES-\Alg$ is $E$--local if and only if it is weakly equivalent $L_EA$, for some $A \in \Alg$.  More precisely, $B$ is $E$--local if and only if the natural weak map $L_E(I(B) \vee S) \xra{\sim} B$ is an equivalence.  Thus $L_E(\Alg)$ is equivalent to the category $L_E^{\prime}(\Alg)$ with objects $A \in \Alg$, and with morphisms from $A$ to $B$ equal to the $L_ES-\Alg$ maps from $L_EA$ to $L_EB$. \\
\end{rem}

\begin{defns}  We define functors with domain $L_E^{\prime}(\Alg)$ as follows.
\begin{enumerate}
\item Let $taq^E(A) = L_E(taq_{L_ES}(L_EA))$.  
\item Define the natural tower of fibrations in $L_E(\Alg)$ under $L_EA$, 
\begin{equation}\label{L_E algebra tower}
\xymatrix{
&&& \vdots \ar[d] \\
&&& P_2^E A\ar[d]^{p_2^E} \\
&&& P_1^E A \ar[d]^{p_1^E} \\
A \ar[rrr]^{e_0^E}  \ar[urrr]^{e_1^E} \ar[uurrr]^{e_2^E} &&& P_0^E A, 
}
\end{equation}
to be the tower obtained by applying $L_E$ to the Andr\'e--Quillen tower (\ref{algebra tower}) $\{P_{L_ES,r}(L_EA)\}$.
\item Let $\displaystyle \widehat{L}_EA = \lim_{r \ra \infty} P_r^E(A).$ \\
\end{enumerate}
\end{defns}

We have the following analogue of \propref{tower prop}.  This proposition is a slight elaboration of  \propref{LE tower prop1} of the introduction. \\

\begin{prop} \label{LE tower prop} Let $f:L_EA \ra L_EB$ be a map in $L_E(\Alg)$.  If $taq^E(f): taq^E(A) \ra taq^E(B)$ is a weak equivalence, so is $P^E_r(f): P^E_r(A) \ra P^E_r(B)$ for all $r$, and thus also
$\widehat{f}: \widehat{L}_EA \ra \widehat{L}_EB$.  Thus, in this case, one gets a factorization by weak $L_ES$--algebra maps
\begin{equation*}
\xymatrix{
L_EA \ar[rr]^{canonical} \ar[dr]_{f}  &  &  \widehat{L}_EA, \\
   & L_EB \ar[ur] &    \\
}
\end{equation*}
where the unlabelled weak map is $L_EB \ra \widehat{L}_EB \xla[\sim]{\widehat{f}} \widehat{L}_EA$.
\end{prop}
\begin{proof} This is proved by induction on $r$, using \lemref{Dr Eiso lemma} and \thmref{tower theorem}(2).
\end{proof}

This is given added punch when combined with the next proposition.\\

\begin{prop} \label{LE tower theorem2} Given $A \in \Alg$ there are natural weak equivalences
\begin{enumerate}
\item $taq^E(A) \simeq L_E(taq(A))$.
\item $\{P_r^E(A)\} \simeq \{L_EP_r(A)\}$, as towers.
\end{enumerate}
\end{prop}
\begin{proof}  First note that if $X$ is an $S$--module, then each of the maps
$$ X \ra L_ES \sm X \ra L_EX$$
is an $E_*$--isomorphism.  

To prove (1), we have
\begin{equation*}
\begin{split}
taq^E(A) & = L_E(taq_{L_ES}(L_EA)) \\
  & \simeq  L_E(taq_{L_ES}(L_ES \sm A))  \text{ by \lemref{taq Eiso lemma}} \\
  & \simeq  L_E(L_ES \sm taq(A))  \text{ by \lemref{taq change of rings}} \\
  & \simeq L_E(taq(A)).
\end{split}
\end{equation*}

To prove (2), we have
\begin{equation*}
\begin{split}
P_r^E(A) & = L_E(P_{L_ES,r}(L_EA)) \\
  & \simeq  L_E(P_{L_ES,r}(L_ES \sm A)) \text{ by \corref{Pr Eiso cor}} \\
  & \simeq  L_E(L_ES \sm P_r(A)) \text{ by \thmref{tower theorem}(3)} \\
  & \simeq L_E(P_r(A)).
\end{split}
\end{equation*}
 
\end{proof}

In constrast to the equivalences in this last proposition, we note that it is {\em not} necessarily true that the natural weak map
$$ L_E(\widehat{A}) \ra \widehat{L}_E(A)$$
is an equivalence (particularly when $E$ is not connective), and thus the convergence of the localized tower is very problematic.  For example, if $A = \PP X$, this map has the form
$$ L_E \left(\prod_{r=0}^{\infty} D_rX\right) \ra \prod_{r=0}^{\infty} L_E D_r X,$$
which would usually not be an equivalence.

\subsection{Proof of the properties of the Andr\'e--Quillen tower} \label{proof of the tower properties}

In the series of papers \cite{goodwillie1, goodwillie2, goodwillie3}, Tom Goodwillie has developed his theory of polynomial resolutions of homotopy functors.  Although \cite{goodwillie3} only explicitly studies such resolutions of functors 
$$ F: \A \ra \B$$
with $\A$ and $\B$ either spaces or spectra, essentially everything in his paper makes sense in a much broader setting.  In particular, his concepts and constructions certainly make sense if $\A$ and $\B$ are (based) topological model categories (a model category tensored over based topological spaces, satisfying properties as in \cite[VII.4]{ekmm}), and $\B$ furthermore is a category in which a directed hocolimit of homotopy cartesion cubical diagrams is again a homotopy cartesion cubical diagram.  One such category is $R-\Alg$.  Another is $R-\Module$, the category of $R$--modules.  

Recall that the tower $\{P_{R,r}(\hspace{.1in})\}$ is defined to be the Goodwillie tower associated to the identity functor on $R-\Alg$. In this subsection we indicate why this tower has the properties given in \thmref{tower theorem}.  We do this by summarizing the main points of Goodwillie's work as they apply to \thmref{tower theorem}.  Throughout we are citing the version of \cite{goodwillie3} of June, 2002.

As in \cite{goodwillie2}, a functor is said to be {\em $r$--excisive} if it takes strongly homotopy cocartesion $(r+1)$--cubical diagrams to homotopy cartesian cubical diagrams. In \cite{goodwillie3}, given a functor $F$, the tower $\{P_rF(\hspace{.1in})\}$ is defined so that $F \ra P_rF$ is the universal arrow to an $r$--excisive functor, up to weak equivalence.   

Goodwillie proves the existence of such a tower by an explicit construction which amounts to modifying $F$ so as to visibly {\em force} certain strongly homotopy cocartesion $(r+1)$--cubical diagrams to transform to homotopy cartesian diagrams.  Readers looking to apply his paper in the setting where the domain and range of $F$ are topological model categories should write `$U_+ \otimes X$' whenever Goodwillie writes `$X \times U$', with $U$ a finite set, and also recall that the domain category has an initial/terminal object.   \cite[Theorem 1.8]{goodwillie3} says that the tower constructed as he describes has the desired universal properties.

For example, there is a strongly cocartesion diagram
\begin{equation*}
\xymatrix{
S^0 \ar[d] \ar[r] &
D^+ \ar[d]  \\
D^- \ar[r] &
S^1,
}
\end{equation*}
representing the circle as the union of two $1$--disks $D^+$ and $D^-$.  Then $P_1F(X)$ is defined to be the homotopy colimit of 
$$ F(X) \ra T_1F(X) \ra T_1T_1F(X) \ra \dots $$
where $T_1F(X)$ is the homotopy pullback of 
\begin{equation*}
\xymatrix{
&
F(D^+ \otimes X) \ar[d]  \\
F(D^- \otimes X) \ar[r] &
F(S^1\otimes X).
}
\end{equation*}
Since $D^+$ and $D^-$ are contractible, $F(D^+ \otimes X)$ and $F(D^- \otimes X)$ are equivalent to the initial/terminal object in the domain category, so that 
$$T_1F(X) \simeq \Omega F(S^1 \otimes X)$$
and 
$$P_1F(X) \simeq \hocolim_{n \ra \infty} \Omega^n F(S^n \otimes X).$$

Already, just using this part of the theory, various parts of \thmref{tower theorem} are evident.  Statement (1), saying that $P_{R,0}A \simeq R$, is clear.  
Statement (3), the change of rings formula, is also clear, noting that $R \sm (\hspace{.1in})$ takes strongly homotopy cocartesion cubes of $S$--algebras to strongly homotopy cocartesion cubes of $R$--algebras, and homotopy cartesion cubes of $S$--algebras to homotopy cartesion cubes of $R$--algebras (see below).  Finally, part of statement (2), that the fiber of $p_1: P_{R,1}A \ra P_{R,0}A$ is naturally equivalent to $taq_R(A)$, follows from the above description of $P_1F$.  

In \cite{goodwillie2, goodwillie3}, Goodwillie develops general theory and examples allowing for connectivity estimates to be made for the maps $F(X) \ra P_rF(X)$ in terms of the connectivity of $X$.  In particular, statement (4), stating that $e_r: A \ra P_{R,r}A$ is $r$--connected if $I(A)$ is 0--connected, can be deduced from \lemref{convergence lemma}.

We are left needing to show the rest of statement (2): that $D_{R,r}(A)$, the fiber of $p_r$, is weakly equivalent to the $r^{th}$ extended power of $taq(A)$, the fiber of $p_1$.  

To show this, we begin by noting that homotopy pullback diagrams in $R-\Alg$ are just diagrams in $R-\Alg$ that are homotopy pullbacks in $R$--modules.  Thus the tower $\{P_{R,r}(\hspace{.1in})\}$, with algebra structures forgotten, is the Goodwillie tower of the inclusion functor
$$ \mathcal I:R-\Alg \hra R-\Module.$$

The category $R-{\text{modules}}$ is a stable model category, in the sense of \cite{hovey2}; in particular, homotopy cocartesian cubical diagrams are equivalent to homotopy cartesion cubical diagrams.  Thus one can apply Goodwillie's analysis in \cite{goodwillie3} of how $D_rF(A)$, the fiber of $P_rF(A) \ra P_{r-1}F(A)$, can be computed by means of cross effects.

The bits of the general theory we need are the following.  Let
$$ F: \A \ra \B$$
be a functor between topological model categories as above, with $\B$ stable.

Let ${\mathbf  r} = \{1,2,\dots, r\}$.  In \cite[\S 3]{goodwillie3}, $\chi_rF$, the $r^{th}$ cross effect of $F$, is defined to the the functor of $r$ variables given as the total homotopy fiber
\begin{equation} \label{cross effect def} 
\chi_rF(A_1, \dots, A_r) = \totfib_{T \subset \mathbf r} F(\coprod_{i \in {\mathbf r} - T} A_i).
\end{equation}

Then \cite[Theorem 6.1]{goodwillie3} says that $D^{(r)}F$, the $r^{th}$ multilinearization of $F$, can be computed by the formula
\begin{equation} \label{stabilization}
 D^{(r)}F(A_1, \dots, A_r) \simeq \hocolim_{n_i \ra \infty} \Omega^{n_1 + \dots + n_r} \chi_rF(S^{n_1} \otimes A_1, \dots, S^{n_r}\otimes A_r).
\end{equation}

Finally, \cite[Theorem 3.5]{goodwillie3} says that there is a natural weak equivalence
\begin{equation} \label{Dr formula}
 D_rF(A) \simeq (D^{(r)}F(A, \dots, A))_{h\Sigma_r}.
\end{equation}

We apply this theory to the case when $F = \mathcal I$.  

Since the coproduct in $R-\Alg$ is just the smash product $\sm_R$, we have
$$ \chi_r\mathcal I(A_1, \dots, A_r) = \totfib_{T \subset \mathbf r} (\bigwedge_{i \in {\mathbf r} - T} A_i). $$

For example, $\chi_2(A,B)$ is the total homotopy fiber of the square

\begin{equation*}
\xymatrix{
A \sm_R B \ar[d] \ar[r] &
A \ar[d]  \\
B \ar[r] &
R
}
\end{equation*}

Recall that the total homotopy fiber is isomorphic to the iterated homotopy fiber.  This makes the next lemma easy to check.

\begin{lem} The natural map
$$ I(A_1) \sm_R \dots \sm_R I(A_r) \ra \chi_r\mathcal I(A_1, \dots, A_r)$$
is an equivalence.
\end{lem}

\begin{cor}  There are natural weak equivalences of $R$--modules

$$ D^{(r)}\mathcal I(A_1, \dots, A_r) \simeq taq(A_1) \sm_R \dots \sm_R taq(A_r),$$
and
$$ D_r\mathcal I(A) = (taq(A)^{\sm_R r})_{h\Sigma_r}.$$
\end{cor}

We have finished our proof of \thmref{tower theorem}, as this last equivalence is just a restatement of the remaining unproven part of statement (2).

\section{Proof of \thmref{classic theorem}, \thmref{susp thm}, \thmref{Tn theorem}, and \thmref{HQ theorem}} \label{proofs section}

In this section we use the theory developed in  the last section to prove the splitting theorems of the introduction.

\subsection{Proof of \thmref{classic theorem}}

\begin{defn}  If $Z$ is a space, let $s(Z): \PP(\Sinfty Z) \ra \Sinfty (QZ)_+$ be the natural weak map in $\Alg$ induced by the weak natural map of $S$--modules
$$ \Sinfty Z \xra{\Sinfty \eta(Z)} \Sinfty QZ \xla{\sim} I(\Sinfty(QZ)_+).$$
\end{defn}

We restate \thmref{classic theorem}.

\begin{thm}  \label{classic theorem1} For all spaces $Z$, the map $s(Z)$ induces an isomorphism on completions, and thus there is a natural factorization of weak algebra maps
\begin{equation*}
\xymatrix{
\PP(\Sinfty Z) \ar[rr]^{canonical} \ar[dr]_{s(Z)}  &  &  \widehat{\PP}(\Sinfty Z). \\
   & \Sinfty (QZ)_+ \ar[ur]_{t(Z)} &    \\
}
\end{equation*}
If $Z$ is 0--connected then all of these maps are weak equivalences.
\end{thm}
\begin{proof}  By \propref{tower prop}, we just need to show that 
$$taq(s(Z)): taq(\PP(\Sinfty Z)) \ra taq(\Sinfty(QZ)_+)$$
is an equivalence.  To see this, consider the diagram:

\begin{equation*}
\xymatrix{
\Sinfty Z \ar[ddd]^{\wr} \ar[ddr]^i \ar[drr]^{\eta} \ar@{=}[rrr] & & & \Sinfty Z  \ar[ddd]^{\wr} \\
&& \Sinfty QZ \ar[ur]^{\epsilon} &  \\
& I(\PP(\Sinfty Z)) \ar[r]^{I(s(Z))} \ar[dl] & I(\Sinfty(QZ)_+) \ar[u]_{\wr} \ar[dr]  &  \\
taq(\PP(\Sinfty Z)) \ar[rrr]^{taq(s(Z))} & & & taq(\Sinfty(QZ)_+) 
}
\end{equation*}

\propref{P prop}(2) says that the left edge is an equivalence, and the left triangle commutes.  Similarly, \propref{SLX prop}(4) shows that the right edge is an equivalence, and the right quadrilateral commutes.  Naturality shows that the bottom quadrilateral commutes.  The map $s(Z)$ was {\em defined} so that the middle quadrilateral commutes.  Finally, the top triangle is just the categorical factorization (\ref{space lifting}).  Thus the diagram commutes, and inspection of the outside square shows that $taq(s(Z))$ can be identified with the identity map on $\Sinfty Z$.
\end{proof}

\begin{rems} \label{classic thm remarks} When $Z$ is connected, the realization that a weak equivalence like $s(Z)$ can be taken to be $E_{\infty}$ dates back to the late 1970's, with the first proof based on a point set analysis of James--Hopf maps \cite{cmt}. Proofs working on the spectrum level were given in \cite[Thm.VII.5.5]{lmms} and  \cite[Prop.4.3]{k1}. Both of these references construct $s(Z)$ using ideas of \cite{rcohen}: see \cite[Appendix B]{k3} for an updated account.

If $Z$ is {\em not} 0--connected, then the fact $t(Z)$ can be constructed to be  $E_{\infty}$ seems to be new, though slightly weaker results were proved in  \cite{cmt}.  We also manage to show the existence of $t(Z)$ without appealing to properties of group completions: again this is new. 

Our proof does {\em not} use the combinatorial approximation $CZ \ra QZ$, though some of the ideas behind that approximation are obviously lurking in the proofs of needed properties of $\Sinfty(QZ)_+$.  We relate our constructions to those using $CZ$ in \appref{splitting appendix}.  
\end{rems}

\subsection{Proof of \thmref{susp thm}}

Suppose that $Z$ is a connected space, and $f: \Sinfty Z \ra X$ a map inducing an isomorphism on $E_*$--homology.  The maps in $\Alg$,
$$ \PP(\Sinfty Z) \xra{s(Z)} \Sinfty(QZ)_+ \xra{\Sinfty \Oinfty f} \Sinfty(\Oinfty X)_+,$$
induce a commutative diagram  
\begin{equation*}
\xymatrix{
L_E \PP(\Sinfty Z) \ar[d] \ar[rr]^-{\sim} && L_E \Sinfty(QZ)_+ \ar[d] \ar[rr]^-{L_E\Sinfty (\Oinfty f)_+} && L_E \Sinfty(\Oinfty X)_+ \ar[d]  \\
\widehat{L}_E \PP(\Sinfty Z) \ar[rr]^-{\sim} && \widehat{L}_E \Sinfty(QZ)_+ \ar[rr]^-{\sim} &&  \widehat{L}_E \Sinfty(\Oinfty X)_+.
}
\end{equation*}
In this diagram, the top left horizontal maps is an equivalence by \thmref{classic theorem}, while the bottom left horizontal map is similarly an equivalence using \propref{LE tower prop1}. The lower right horizontal map is an equivalence by \propref{LE tower prop1}, since our hypothesis that $f$ is an $E_*$--isomorphism implies that $taq^E(f)$ is an equivalence.

As the left vertical map is an $E_*$--monomorphism, we conclude that so is the upper right horizontal map.  Otherwise said,
$$ (\Oinfty f)_*: E_*(QZ) \ra E_*(\Oinfty X)$$
is monic. 

\subsection{Telescopic functors, and the proof of \thmref{Tn theorem}}

We fix a prime $p$ and work $p$--locally. For $n \geq 1$, let $K(n)$ be the $n^{th}$  Morava K--theory, and $T(n)$ be the telescope of a $v_n$--self map of a finite complex of type $n$.  It is known that the Bousfield class of $T(n)$ is independent of choice of $v_n$--self map, and that $\langle K(n) \rangle \leq \langle T(n) \rangle$, with equality holding if and only if the Telescope Conjecture holds. (See \cite{bousfield4} for background and more references.)

\begin{thm}  There exists a functor $\Phi_n: \text{Spaces} \ra \text{$S$--modules}$ and a natural weak equivalence $\Phi_n \Oinfty X \simeq L_{T(n)}X$.
\end{thm}

With the result stated at the level of homotopy categories, and with $K(n)$ replacing $T(n)$, this is the main theorem of \cite{k2}.  However the sorts of constructions given there, and in \cite{bousfield1} (for $n=1$), yield the theorem as stated.  In particular, in the recent paper \cite{bousfield4}, Bousfield proves the theorem as stated using the model category of spectra of \cite{bf}.  But this category is known \cite{schwedeshipley} to be Quillen equivalent to the $S$--modules of \cite{ekmm}. 

As a corollary, we obtain a proof of (\ref{Tn lifting}), which we restate here.

\begin{cor} \label{Tn lifting cor} (Compare with \cite{k2}.)
There is a natural factorization by weak $S$--module maps
\begin{equation*} 
\xymatrix{   & L_{T(n)}\Sinfty \Oinfty X \ar[dr]^{L_{T(n)}\epsilon(X)} &    \\
L_{T(n)}X \ar@{=}[rr] \ar[ur]^{\eta_n(X)}  &  &  L_{T(n)}X. \\
}
\end{equation*}
\end{cor}
\begin{proof}  Apply $\Phi_n$ to the commutative diagram 
\begin{equation} 
\xymatrix{   & \Oinfty \Sinfty \Oinfty X \ar[dr]^{\Oinfty\epsilon(X)} &    \\
\Oinfty X \ar@{=}[rr] \ar[ur]^{\eta(\Oinfty X)}  &  &  \Oinfty X. \\
}
\end{equation}
\end{proof}

\begin{defn}  If $X$ is an $S$--module, let $$s_n(X): L_{T(n)}\PP(X) \ra L_{T(n)}\Sinfty (\Oinfty X)_+$$ be the natural weak map in $L_{T(n)}(\Alg)$ induced by the weak natural map of $S$--modules
$$ L_{T(n)}X \xra{\eta_n(X)} L_{T(n)}\Sinfty \Oinfty X \xla{\sim} L_{T(n)}I(\Sinfty(\Oinfty X)_+).$$
\end{defn}

 We restate \thmref{Tn theorem}.

\begin{thm}  For all spectra $X$, the map $s_n(X)$ induces an isomorphism on $L_{T(n)}$--completions, and thus there is a natural factorization of weak algebra maps
\begin{equation*}
\xymatrix{
L_{T(n)}\PP(X) \ar[rr]^{canonical} \ar[dr]_{s_n(X)}  &  &  \widehat{L}_{T(n)}\PP(X). \\
   & L_{T(n)}\Sinfty (\Oinfty X)_+ \ar[ur]_{t_n(Z)} &    \\
}
\end{equation*}
\end{thm}
\begin{proof}  Denote $L_{T(n)}$ by $L$.  By \propref{LE tower prop}, we just need to show that 
$$taq(s_n(X)): Ltaq(\PP(X)) \ra Ltaq(\Sinfty(\Oinfty X)_+)$$
is an equivalence.  To see this, consider the diagram:
{\small 
\begin{equation*}
\xymatrix{
LX \ar[ddd]^{\wr} \ar[ddr]^i \ar[drr]^{\eta_n} \ar@{=}[rrr]  & & & LX   \ar[ddd]^{\wr} \\
&& L\Sinfty \Oinfty X \ar[ur]^{L\epsilon} &  \\
& LI(\PP(X)) \ar[r]^-{I(s_n(X))} \ar[dl]  & LI(\Sinfty(\Oinfty X)_+) \ar[u]_{\wr} \ar[dr]  &  \\
Ltaq(\PP(X)) \ar[rrr]^{taq(s_n(X))}  & & & Ltaq(\Sinfty(\Oinfty)_+) 
}
\end{equation*}}

The top triangle commutes by \corref{Tn lifting cor}, and the map $s_n(X)$ was defined so that the middle quadrilateral commutes.  The rest of the diagram commutes for the same reasons as in the proof of \thmref{classic theorem}.
\end{proof}

\subsection{Rational spectra and the proof of \thmref{HQ theorem}}

We restate (\ref{HQ lifting}) as a lemma.

\begin{lem} \label{HQ lifting lemma} For all $0$--connected spectra $X$, there is a natural factorization by weak $S$--module maps
\begin{equation*} 
\xymatrix{   & L_{H\Q}\Sinfty \Oinfty X \ar[dr]^{L_{H\Q}\epsilon(X)} &    \\
L_{H\Q}X \ar@{=}[rr] \ar[ur]^{\eta_0(X)}  &  &  L_{H\Q}X. \\
}
\end{equation*}
For $X$ just $-1$--connected, such a factorization also exists, but can not be made to be natural. 
\end{lem}

Assuming this for the moment, we continue as we did in the last subsection.

\begin{defn}  If $X$ is a $-1$--connected $S$--module, let $$s_0(X): L_{H\Q}\PP(X) \ra L_{H\Q}\Sinfty (\Oinfty X)_+$$ be the weak map in $L_{H\Q}(\Alg)$ induced by the weak natural map of $S$--modules
$$ L_{H\Q}X \xra{\eta_0(X)} L_{H\Q}\Sinfty \Oinfty X \xla{\sim} L_{H\Q}I(\Sinfty(\Oinfty X)_+).$$
\end{defn}

 We restate \thmref{HQ theorem}.

\begin{thm}  For all $-1$--connected spectra $X$, the map $s_0(X)$ induces an isomorphism on $L_{H\Q}$--completions, and thus there is a factorization of weak algebra maps
\begin{equation*}
\xymatrix{
L_{H\Q}\PP(X) \ar[rr]^{canonical} \ar[dr]_{s_0(X)}  &  &  \widehat{L}_{H\Q}\PP(X). \\
   & L_{H\Q}\Sinfty (\Oinfty X)_+ \ar[ur]_{t_0(Z)} &    \\
}
\end{equation*}
This is natural for $0$--connected $X$, and, in this case, all three maps are equivalences.
\end{thm}

The theorem follows from the lemma in the usual way.

\begin{proof}[Proof of \lemref{HQ lifting lemma}]  The idea behind this lemma is that passing to homotopy groups is a full and faithful process on the homotopy category of $H\Q$--local $S$--modules.  

A version of the lemma then follows easily as follows. Let 
$$ \eta(X)_*: \pi_*(X) \ra \pi_*(\Sinfty \Oinfty X) $$
be the map induced by the canonical inclusion
$$ \eta(X): \Oinfty X \ra Q \Oinfty X.$$
Note that $\eta(X)_*$ is a map of abelian groups if $*>0$, but is only a map of sets if $*=0$.  However, for all $* \geq 0$, we have that $\epsilon(X)_* \circ \eta_*(X)$ is the identity.

Thus if $X$ is $0$--connected, there is a canonical natural homotopy class of maps
$$\eta_0(X): L_{H\Q} X \ra L_{H\Q}\Sinfty \Oinfty X$$
realizing $\eta(X)_* \otimes \Q$.  

If $X$ is just $-1$--connected, one can still choose a $\Q$--linear section to $\epsilon(X)_* \otimes \Q$, and then realize this, defining $\eta_0(X)$.  However, this cannot possibly be taken to be natural by the following argument, which the author learned from Pete Bousfield. If $\eta_0(X)$ were natural, then the maps $\eta_0(HV)$, with $V$ a $\Q$--vector space, would define a natural section to the augmentation
$$ \Q[V] \ra V$$
defined on the category of $\Q$--vector spaces.  (Here $\Q[V]$ denotes the $\Q$--vector space with basis $V$.)  But it is well known, and easily verified, that there exist no nontrivial natural transformations $V \ra \Q[V]$.

A careful reader may be wondering if, in the 0--connected case, one can construct a natural weak section at the model category level, and not just a natural section in the homotopy category.  This is also possible: one can apply \cite[Theorem 4]{mccarthy1}, which implies that if $P^1(X)$ is the codegree 1 approximation to $L_{H\Q}(\Sinfty \Oinfty X)$ (in the sense of dual calculus), then the composite 
$$ P^1(X) \ra L_{H\Q} \Sinfty \Oinfty X \xra{L_{H\Q}\epsilon(X)} L_{H\Q} X$$
is an equivalence for $0$--connected $X$.
\end{proof}

\section{When $s_n(X)$ is an equivalence, and related matters} \label{snX section}

Recall that $\Sp_n = \{ \text{$S$--modules } X \ | \ s_n(X) \text{ is an equivalence} \}$ and that  $\Sp^K_n = \{ \text{$S$--modules } X \ | \ s^K_n(X) \text{ is an equivalence} \}$.   Recall also that $c(n)$ denotes the smallest integer $c$ such that $\tilde{T}(n)_*(K(\mathbb Z/p,c)) = 0$.  The analogous integer associated to $K(n)_*$ is $n+1$: by the calculations in \cite{rw}, $\tilde{K}(n)_*(K(\mathbb Z/p,c)) = 0$ if and only if $c \geq n+1$.\footnote{When $p=2$, the reference is \cite[Appendix]{johnson wilson}.}

The starting point for the detailed results about $\Sp_n$ and $\Sp_n^K$ given in \secref{sn iso subsection} is the following partial result. \\
 
\begin{thm} \label{sn(X) iso thm 2}  Let $X$ be an $S$--module such that
$L_{n-1}^f X \simeq *$.

\noindent (1) $X \in \Sp_n$ if $X$ is $c(n)$--connected. 

\noindent (2) $X \in \Sp_n^K$ if $X$ is $n+1$--connected. \\
\end{thm}

The proof of this theorem is slightly long and delicate.  We organize it into the following steps: \\

\noindent{\bf Step 1} \ We show that if $F$ is a finite $S$--module of type $n$, then $\Sigma^d F \in \Sp_n$ for $d \gg 0$. \\

\noindent{\bf Step 2} \ Assuming Step 1, we show that if $F$ is a 0--connected finite $S$--module of type $n$, then $\Sigma^{c(n)} F \in \Sp_n$, and $\Sigma^{n+1} F \in \Sp_n^K$. \\

\noindent{\bf Step 3} \  We show that the classes $\Sp_n$ and $\Sp_n^K$ are closed under various constructions. \\

\noindent{\bf Step 4} \  We show that, starting from the finite $S$--modules shown in Step 2 to be in $\Sp_n$ and $\Sp_n^K$, one can build all $X$'s as in \thmref{sn(X) iso thm 2} using the constructions of Step 3. \\

These will be proved in the next four subsections.  Then, armed with \thmref{sn(X) iso thm 2}, we will systematically work our way through the proofs of the various results stated in subsections \secref{sn iso subsection} and \secref{more homo cor subsection}.

\subsection{Proof of \thmref{sn(X) iso thm 2}: step 1}

\begin{lem} \label{s = sn lem}   Let $Z$ be a space whose suspension spectrum is finite of type at least $n$.  Then, for $d \gg 0$, $L_{T(n)}s(\Sigma^d Z) \simeq s_n(\Sinfty \Sigma^dZ)$.
\end{lem}

Postponing the proof momentarily, we note the following corollary.

\begin{cor} \label{big F in Sn lem}  Let $F$ be a finite of type at least $n$.  Then for $d \gg 0$, $s_n(\Sigma^d F)$ is an equivalence.
\end{cor}

\begin{proof}  Replacing $F$ by a high suspension of $F$ if needed, we can assume that $F \simeq \Sinfty Z$, for some space $Z$.  By the lemma, for large enough $d$, $s_n(\Sigma^d F)$ will be homotopic to $L_{T(n)}s(\Sigma^d Z)$, which is an equivalence by \thmref{classic theorem}.
\end{proof}

\begin{proof}[Proof of \lemref{s = sn lem}]  Let $Z$ be a space.  By construction, the two maps of $L_{T(n)}S$--algebras
$$ s_n(\Sinfty Z), L_{T(n)}s(Z): L_{T(n)}\PP(\Sinfty Z) \ra L_{T(n)}\Sinfty(QZ)_+$$
will be homotopic if and only if the two maps of $L_{T(n)}S$--modules
$$ \eta_n(\Sinfty Z), L_{T(n)}\Sinfty \eta(Z): L_{T(n)}\Sinfty Z \ra L_{T(n)}\Sinfty QZ$$
are homotopic.

Now suppose that $\Sinfty Z$ has type at least $n$, and is at least $0$--connected.  Then $\Sinfty Z$ is $T(i)_*$--acyclic for $0 \leq i \leq n-1$, and thus the same is true for $\Sinfty QZ$, since by \thmref{classic theorem}, $\Sinfty (QZ)_+ \simeq \PP(\Sinfty Z)$.

If a spectrum $X$ is $T(i)_*$--acyclic for $0 \leq i \leq n-1$, then $L_{T(n)}X \simeq L_n^f X \simeq L_n^f S \sm X$.  Applying this to our situation, we see that  $\eta_n(\Sinfty \Sigma^dZ)$ and $L_{T(n)}\Sinfty \eta(\Sigma^dZ)$ correspond to homotopy classes of maps of $S$--modules  
$$ x_n(d), x(d) \in [\Sinfty Z ,  L_n^f S \sm \Sigma^{-d} \Sinfty Q\Sigma^d Z].$$

We now show that, if $d$ is very large, then $x(d) = x_n(d)$. This follows from three observations.  

Firstly, the naturality of $\eta$ and $\eta_n$ ensures that $x(d)$ maps to $x(d+1)$, and $x_n(d)$ maps to $x_n(d+1)$, under the homomorphism
$$ [\Sinfty Z ,  L_n^f S \sm \Sigma^{-d} \Sinfty Q\Sigma^d Z] \ra [\Sinfty Z ,  L_n^f S \sm \Sigma^{-(d+1)} \Sinfty Q\Sigma^{d+1} Z]
$$
induced by the evaluation map $\Sigma Q\Sigma^dZ \ra Q\Sigma^{d+1}Z$.  

Secondly, the colimit 
$$\colim_{d \ra \infty}[\Sinfty Z ,  L_n^f S \sm \Sigma^{-d} \Sinfty Q\Sigma^d Z]$$
can be identified with $[\Sinfty Z ,  L_n^f S \sm \Sinfty  Z]$.

Thirdly, the key properties of $\eta$ and $\eta_n$, (\ref{space lifting}) and (\ref{Tn lifting}), imply that under this identification, $\displaystyle \colim_{d \ra \infty} x(d)$ and $\displaystyle \colim_{d \ra \infty} x_n(d)$
each correspond to the canonical element: the unit $S \ra L_n^f S$ smashed with the identity on $\Sinfty Z$.

We conclude that $\displaystyle \colim_{d \ra \infty}(x(d)-x_n(d))$ is zero in the colimit, and thus $x(d)-x_n(d)$ is zero at a finite stage of this colimit.  Otherwise said, for $d$ large, we have $x(d) = x_n(d)$.
\end{proof}

\subsection{Proof of \thmref{sn(X) iso thm 2}: step 2}

\begin{prop} \label{type n+1 is strongly acyclic prop}  Let $F$ be a $-1$--connected finite $S$--module of type $n+1$.

\noindent (1) $\Tilde{T}(n)_*(\Oinfty \Sigma^{c(n)}F) = 0$.

\noindent (2) $\Tilde{K}(n)_*(\Oinfty \Sigma^{n+1}F) = 0$.
\end{prop}

As a corollary, one deduces the assertion of step 2. 

\begin{cor} \label{finites in Sn prop} Let $F$ be a 0--connected finite $S$--module of type $n$.  

\noindent (1) $\Sigma^{c(n)}F \in \Sp_n$.

\noindent (2) $\Sigma^{n+1}F \in \Sp_n^K$.
\end{cor}

\begin{proof}  Let $F$ be 0--connected and finite of type $n$.  Let $v: \Sigma^d F \ra F$ be a $v_n$--map, an isomorphism in both $T(n)_*$ and $K(n)_*$.  Let $W$ be the fiber.  By statement (1) of the proposition, $\Oinfty \Sigma^{c(n)}W$ is $T(n)_*$--acyclic.  By an Atiyah--Hirzebruch--Serre spectral sequence argument, it follows that $\Sigma^c v$ induces an isomorphism
$$ v_*: T(n)_*(\Oinfty \Sigma^{d+c} F) \xra{\sim} T(n)_*(\Oinfty \Sigma^c F)$$
if $c \geq c(n)$.  

Now consider the diagram 
\begin{equation*}
\xymatrix{
L\PP(\Sigma^{dN+c(n)}F) \ar[d] \ar[rr]^-{s_n(\Sigma^{dN+c(n)}F)} & & L\Sinfty (\Oinfty(\Sigma^{dN+c(n)}F))_+ \ar[d]  \\
L\PP(\Sigma^{c(n)}F) \ar[rr]^-{s_n(\Sigma^{c(n)}F)} & & L\Sinfty(\Oinfty \Sigma^{c(n)}F)_+
}
\end{equation*}
where $L = L_{T(n)}$, and the vertical maps are induced by $v^N$.  The above discussion implies that the right vertical map is an equivalence.  The left vertical map is an equivalence because $v$ is a $T(n)_*$--isomorphism.  Finally, for large enough $N$ the top horizontal map is an equivalence as a consequence of Step 1. Thus the bottom map is an equivalence, proving statement (1) of the corollary.  Statement (2) is proved similarly.
\end{proof}

To prove \propref{type n+1 is strongly acyclic prop}, we use ideas from \cite{hrw}. 

If $H\Z/p$ is strongly $E_*$--acyclic, let $c_p(E)$ be the smallest $c$ such that $\Tilde{E}_*(K(\Z/p, c)) = 0$. The two statements of \propref{type n+1 is strongly acyclic prop} are then just special cases of the following proposition. \\

\begin{prop}  \label{strongly acyclic prop} Suppose a $-1$--connected spectrum $X$ has $p$--torsion homotopy groups.  If $X$ is strongly $E_*$--acyclic, then $\Tilde{E}_*(\Oinfty \Sigma^{c_p(E)} X) = 0$.
\end{prop}

\begin{proof}  We argue as in \cite[\S 3]{hrw}.  

Firstly, the argument proving \cite[Prop.3.4]{hrw} shows that if $\Oinfty \Sigma^{c} X$ is $E_*$--acyclic, then so is $\Oinfty \Sigma^{c} (X \sm Y)$ for any $-1$--connected spectrum $Y$.  

Now let statement (c) be the statement that $\Tilde{E}_*(\Oinfty \Sigma^{c} X) = 0$.  By assumption, statement (c) is true if $c$ is very large.  We complete the proof of the proposition by showing that, if $c \geq c_p(E)$, then statement (c+1) implies statement (c).

Consider the fibration sequence of $S$--modules
$$ \bar{H} \ra S \ra H\Z.$$
This induces a fibration sequence of spaces
$$ \Oinfty \Sigma^{c+1} (X \sm \Sigma^{-1}\bar{H}) \ra \Oinfty \Sigma^c X \ra \Oinfty \Sigma^c (X \sm H\Z).$$
By our inductive assumption, and the observation above, the first of these spaces is $E_*$--acyclic.  Thus the second map is an $E_*$--isomorphism.  But $\Oinfty \Sigma^c (X \sm H\Z)$ will be a weak product of Eilenberg--MacLane spaces of type $K(A,c)$ where $A$ is $p$--torsion and $c \geq c_p(n)$.  Thus this space is also $E_*$--acyclic, and we conclude the same for $\Oinfty \Sigma^c X$.
\end{proof}

\subsection{Proof of \thmref{sn(X) iso thm 2}: step 3}  The following proposition says that $\Sp_n$ and $\Sp_n^K$ are closed under various constructions. \\

\begin{prop} \label{Sn prop}  Let $\Sp^{\prime}$ be either $\Sp_n$ or $\Sp_n^K$. 

\noindent (1)  $X \in \Sp^{\prime}$ if and only if $X \langle -1 \rangle \in \Sp^{\prime}$.

\noindent (2) If $X \in \Sp^{\prime}$, then $\pi_0(X) = 0$.

\noindent (3)  Let $X \ra Y \ra Z$ be a cofibration sequence of $S$--modules with $X$ and $Y$ $-1$--connected and $Z$ $0$--connected.  Then $X,Y \in\Sp^{\prime}$ implies that $Z \in \Sp^{\prime}$.  In particular, if $X \in \Sp^{\prime}$ is $-1$--connected, then $\Sigma X \in \Sp^{\prime}$.

\noindent (4)  Let $X$ be the filtered homotopy colimit of $S$--modules $X_i$.  If $X_i \in \Sp^{\prime}$ for all $i$, then $X \in \Sp^{\prime}$. \\
\end{prop}

\begin{proof}  We will prove the statements when $\Sp^{\prime} = \Sp_n$; the proofs when $\Sp^{\prime} = \Sp_n^K$ are similar.  Let $L$ denote $L_{T(n)}$.

The first two properties follow from the fact that Eilenberg--MacLane spectra are $T(n)_*$--acyclic.  

In more detail, $X\langle -1 \rangle \ra X$ is a $T(n)_*$--isomorphism, and thus so is $ \PP(X\langle -1 \rangle) \ra \PP(X)$.  Also $\Oinfty X\langle -1 \rangle = \Oinfty X$.  Thus in the square
\begin{equation*}
\xymatrix{
T(n)_*(\PP(X\langle -1 \rangle)) \ar[d] \ar[rr]^-{s_{n*}(X \langle -1 \rangle)} & & T(n)_*(\Oinfty(X\langle -1 \rangle) \ar[d]  \\
T(n)_*(\PP(X)) \ar[rr]^-{s_{n*}(X)} & & T(n)_*(\Oinfty X)
}
\end{equation*}
both vertical maps are equivalences and (1) follows.

For (2), consider the commutative square
the square
\begin{equation*}
\xymatrix{
T(n)_*(\PP (X\langle 0 \rangle)) \ar[d] \ar[rr]^-{s_{n*}(X \langle 0 \rangle)} & & T(n)_*(\Oinfty(X\langle 0 \rangle)) \ar[d]  \\
T(n)_*(\PP(X)) \ar[rr]^-{s_{n*}(X)} & & T(n)_*(\Oinfty X).
}
\end{equation*}
The left vertical map is an isomorphism, and the horizontal maps are monic by \thmref{Tn theorem}.  As $\Oinfty(X\langle 0 \rangle)$ is just one of the path components of $\Oinfty X$, the right vertical map is only epic if $\pi_0(X) = 0$, and (2) follows.

Properties (3) and (4) follow by combining \propref{P prop}, \propref{SLX prop}, and \propref{LE model cat prop}. \\ 
\end{proof}

\subsection{Proof of \thmref{sn(X) iso thm 2}: step 4}

It is convenient to make the following definition, a variant on similar notions in the literature. \\

\begin{defn}  Let $\C$ be any collection of $S$--modules.  Say that an $S$--module $X$ is {\em built from $\C$} if $X \simeq \hocolim_i X_i$, for some sequence $X_0 \ra X_1 \ra \dots$ such that $X_0$ is a wedge of $S$--modules in $\C$, and, for all $i \geq 0$, $X_{i+1}$ is the cofiber of a map $W_i \ra X_i$ with $W_i$ a wedge of $S$--modules in $\C$.  \\
\end{defn}

Note that properties (3) and (4) of \propref{Sn prop} imply that if $\C$ is any subset of $0$--connected $S$--modules in $\Sp_n$ or $\Sp_n^K$, then any $S$--module built from $\C$ is also in $\Sp_n$ or $\Sp_n^K$.

The following proposition is a variant of a well known consequence (as in \cite{miller}) of the Nilpotence Theorems. I would like to thank Pete Bousfield for suggesting the simple proof. \\

\begin{prop} \label{Ln acyclic prop} $L_{n-1}^f X \simeq *$ and $X$ is $c$--connected if and only if $X$ can be built from $c$--connected finite $S$--modules of type $n$. \\
\end{prop}

To prove this we first need a lemma. \\

\begin{lem}  Suppose $L_{n-1}^f X \simeq *$.  Then any $f: Y \ra X$, with $Y$ finite of type at most $n$ and $c$--connected, can be factored as a composite $Y \ra F \ra X$ such that $F$ is a finite $S$--module that is both $c$--connected and of type $n$. 
\end{lem}

\begin{proof}  We prove this by downwards induction on the type of $Y$.  The induction is begun by noting that there is nothing to prove if $Y$ has type $n$.  

So suppose the lemma has been established for all $g: Z \ra X$, where $Z$ is $c$--connected of type at least i+1.  Let $f: Y \ra X$, where $Y$ has type $i$ and is $c$--connected.  By \cite{hs}, there exists a $v_i$--self map $v: \Sigma^d Y \ra Y$.  Since $T(i)_*(X) = 0$, there exists an $N$ such that $\Sigma^{Nd} Y \xra{v^N} Y \xra{f} X$ is null.  Letting $Z$ be the cofiber of $v^N$, it follows that $f$ factors as a composite $Y \ra Z \xra{g} X$.  Since $Z$ is of type $i+1$ and is still $c$--connected, $g$, and thus $f$, factors as needed. \\
\end{proof}

\begin{proof}[Proof of \propref{Ln acyclic prop}]  The Nilpotence Theorem implies that $T(i)_*(F) = 0$ whenever $F$ is finite of type greater than $i$.  Thus one implication is clear:  if $X$ can be built from $c$--connected finite $S$--modules of type at least $n$, then $L_{n-1}^f X \simeq *$ and $X$ is $c$--connected. 

Conversely, suppose that $L_{n-1}^f X \simeq *$ and $X$ is $c$--connected.  We describe how to construct a diagram 
\begin{equation*}
\xymatrix{
X_0 \ar[d]_-{g_0} \ar[r]^{j_0} & X_1 \ar[dl]_-{g_1} \ar[r]^{j_1} & X_2 \ar[r]^{j_2} \ar[dll]_-{g_2} & \dots \\
X & & & 
}
\end{equation*}
showing that $X$ is built from $c$--connected finites of type at least $n$.

First choose a wedge of $c$--connected spheres $T$ and a map $f: T \ra X$ that is onto in $\pi_*$.  By the last lemma, this factors as a composite $T \ra X_0 \xra{g_0} X$ with $X_0$ a wedge of $c$--connected finites of type $n$.

Assume $g_i: X_i \ra X$ has been constructed with $\pi_*(g_i)$ onto, and $X_i$ $c$--connected.  Let $Y_i$ be the fiber of $g_i$. Choose a wedge of $c$--connected spheres $T$ and a map $f: T \ra Y_i$ that is onto in $\pi_*$.  By the last lemma, this factors as a composite $T \ra W_i \xra{g_0} Y_i$ with $W_i$ a wedge of $c$--connected finites of type $n$.  If we then let $X_{i+1}$ be the cofiber of the composite $W_i \ra Y_i \ra X_i$, it follows that $g_i$ will factor as a composite $X_i \xra{j_i} X_{i+1} \xra{g_{i+1}} X$.  

By construction, $\pi_*(g_i)$ is onto and $ker(\pi_*(j_i)) = ker(\pi_*(g_i))$ for all $i$.  It follows that $\hocolim_i X_i \simeq X$ as needed.  \\
\end{proof}

\subsection{Virtual homology and the proof of \lemref{strongly acyclic lemma}}

We need a variant of \propref{strongly acyclic prop}.

If Eilenberg MacLane spectra are strongly $E_*$--acyclic, let $c(E)$ be the smallest $c$ such that $\Tilde{E}_*(K(\Z, c)) = 0$. 

\begin{prop}  \label{strongly acyclic prop 2} If a $-1$--connected spectrum $X$ is strongly $E_*$--acyclic, then $\Tilde{E}_*(\Oinfty \Sigma^{c(E)} X) = 0$.
\end{prop}

The proof of \propref{strongly acyclic prop} goes through without change.

\begin{proof}[Proof of \lemref{strongly acyclic lemma}]

Suppose $X$ is $-1$--connected.  We need to show that 
$ \tilde{E}_*(\Oinfty \Sigma^c X) = 0 $
for large $c$ if and only if
$ \tilde{E}_*(\Oinfty X\langle d \rangle) = 0 $
for large $d$.

Let $P^dX$ denote the $d^{th}$ Postnikov section of $X$, so there is a cofibration sequence of spectra
$$ X\langle d \rangle \ra X \ra P^dX.$$
Then, for all $c \geq 1$, there is a fibration sequence of spaces
$$ \Oinfty \Sigma^{c-1} P^dX \ra \Oinfty \Sigma^c (X \langle d \rangle) \ra \Oinfty \Sigma^c X.$$
If $c > c(E)$, this fiber will be $E_*$--acyclic, and thus there will be an isomorphism
$$ E_*(\Oinfty \Sigma^c (X \langle d \rangle)) \xra{\sim} E_*(\Oinfty \Sigma^c X).$$

Suppose $\tilde{E}_*(\Oinfty X\langle d \rangle) = 0$.  Then $\tilde{E}_*(\Oinfty \Sigma^c (X \langle d \rangle)) = 0$ for all $c \geq 0$.  Thus by our remarks above, $ \tilde{E}_*(\Oinfty \Sigma^c X) = 0$
for $c > c(E)$.

Conversely, suppose $ \tilde{E}_*(\Oinfty \Sigma^c X) = 0$ for all large $c$.  Then, for all $d \geq 0$, 
$\tilde{E}_*(\Oinfty \Sigma^c (X \langle d \rangle)) = 0$ for all large $c$, i.e. $X\langle d \rangle$ is strongly $E_*$--acyclic.  If $d \geq c(E) -1$, we can apply \propref{strongly acyclic prop 2} to the spectrum $\Sigma^{-c(E)}(X\langle d \rangle)$ to conclude that $ \tilde{E}_*(\Oinfty X\langle d \rangle) = 0 $.
\end{proof}

\subsection{Proof of \thmref{sn bar iso thm} and \propref{Lnf strongly acyclic prop}}

We start with a general lemma.

\begin{lem} Suppose $f: X \ra Y$ is a map between $0$--connected spectra with cofiber $C$.  For any homology theory $E_*$, there are implications $(1) \Rightarrow (2) \Rightarrow (3)$.

\noindent (1) $\tilde{E}(\Oinfty \Sigma^{-1}C) = 0$.

\noindent (2) $E_*(\Oinfty f)$ is an isomorphism.

\noindent (3) $\tilde{E}(\Oinfty C) = 0$.
\end{lem}

\begin{proof}  To see that (1) implies (2), note that $\Oinfty \Sigma^{-1}C$ is the fiber of $\Oinfty f$.  

To see that (2) implies (3), note that
$$ \Sinfty(\Oinfty X)_+ \xra{\Oinfty f} \Sinfty(\Oinfty Y)_+ \ra
\Sinfty(\Oinfty C)_+$$
is a cofibration sequence in $\Alg$. Thus \propref{LE model cat prop} applies to say that if $\Oinfty f$ is an $E_*$--isomorphism, then $\Sinfty(\Oinfty C)_+$ is $E_*$--equivalent to $S$.
\end{proof}

\begin{proof}[Proof of \thmref{sn bar iso thm}]  We prove statement (1); the proof of statement (2) is similar.

We temporarily introduce a new class of spectra: let 
$$\tilde{\Sp}_n 
= \{ X \in \Sp \ | \ \Sigma^c(X\langle -1 \rangle) \in \Sp_n \text{ for large } c \}.$$

Let $C_{n-1,d}(X)$ and $f_d$ be defined so that 
$$ C_{n-1,d}(X) \xra{f_d} X\langle d \rangle \ra (L_{n-1}^f X)\langle d \rangle$$
is a fibration sequence of $S$--modules.

As $L_{n-1}^f$ is smashing, and $T(i)_*(T(n)) = 0$ for $0 \leq i \leq n-1$, it follows that $L_{n-1}^f X$ is always $T(n)_*$--acyclic, thus $f_d$ is always a $T(n)_*$--isomorphism.

\thmref{sn(X) iso thm 2} then implies that, for all large $d$ and large $c$, $$C_{n-1,d}(X), \Sigma^c C_{n-1,-1}(X) \in \Sp_n.$$  

Now consider the diagram
\begin{equation*}
\xymatrix{
T(n)_*(\PP(C_{n-1,d}(X))) \ar[d] \ar[r] &
T(n)_*(\Oinfty C_{n-1,d}(X)) \ar[d]^{T(n)_*(\Oinfty f_d)}  \\
T(n)_*(\PP(X \langle d \rangle)) \ar[r] &
T(n)_*(\Oinfty X \langle d \rangle)
}
\end{equation*}
If $d$ is large, both the top map and the left map are isomorphisms, and we conclude that $X\langle d \rangle \in \Sp_n$ if and only if $T(n)_*(\Oinfty f_d)$ is an isomorphism.  Thus $X \in \bar{\Sp}_n$ if and only if $T(n)_*(\Oinfty f_d)$ is an isomorphism for large $d$.

Similarly $X \in \tilde{\Sp}_n$ if and only if $T(n)_*(\Oinfty \Sigma^cf_{-1})$ is an isomorphism for large $c$.

The last lemma and \lemref{strongly acyclic lemma} now combine to say that 
$$X \in \bar{\Sp}_n \Leftrightarrow X \in \tilde{\Sp}_n \Leftrightarrow L_{n-1}^fX \text{ is strongly $T(n)_*$--acyclic.}$$

The inclusion $\Sp_n \subset \tilde{\Sp}_n$ is evident, using \propref{Sn prop}, thus  $\Sp_n \subset \bar{\Sp}_n$.
\end{proof}

\begin{proof}[Proof of \propref{Lnf strongly acyclic prop}]
Suppose condition (1) holds: $T(i)_*(X) = 0$ for $1 \leq i \leq n-1$.  Then $L_{n-1}^f X \ra L_{\Q}X$ is an equivalence.  Since rational spectra are certainly strongly $T(n)_*$--acyclic, condition (2) holds: $L_{n-1}^f X$ is strongly $T(n)_*$--acyclic.

Condition (2) obviously implies condition (3).

Suppose condition (3) holds: $L_{n-1}^f X$ is strongly $K(n)_*$--acyclic.  The main theorem of either \cite{bousfield5} and \cite{wilson} implies that if a spectrum $Y$ is strongly $K(n)_*$--acyclic, then it is strongly $K(i)_*$--acyclic for $1\leq i \leq n$.  Applied to our situation, we deduce that $L_{n-1}^f X$ is strongly $K(n-1)_*$--acyclic.  With notation as in the last proof, \thmref{sn(X) iso thm 2} implies that $C_{n-1,-1}(X)$ is also strongly $K(n-1)_*$--acyclic.  Thus so is $X$, i.e. condition (4) holds.

Finally, if $X$ is strongly $K(n-1)_*$--acyclic, then the Bousfield--Wilson theorem implies that $X$ is $K(i)_*$--acyclic for $1 \leq i \leq n-1$, and thus condition (5) holds.
\end{proof}

\subsection{Proof of \thmref{sn* ses thm} and \thmref{ses thm}}

The fact that the Kunneth Theorem holds for $K(n)_*$ allows for special calculational techniques.  For example, \cite[Thm.7.7]{ekmm} applies to show 
that, if $A \ra B \ra C$ is a cofibration sequence in $\Alg$, the bar spectral sequence converging to $K(n)_*(C)$ has 
$$ E^2_{p,q} = \Tor_{p,q}^{K(n)_*(A)}(K(n)_*(B), K(n)_*).$$

This has the following consequence of relevence to us.

\begin{lem}  Suppose $f: X \ra Y$ is a map between $0$--connected spectra with cofiber $C$.  If $K(n)_*(\Oinfty f)$ is monic, there is a short exact sequence of $K(n)_*$--Hopf algebras
$$ K(n)_*(\Oinfty X) \xra{(\Oinfty f)_*} K(n)_*(\Oinfty Y) \ra K(n)_*(\Oinfty C).$$
\end{lem}
\begin{proof} $K(n)_*(\Oinfty X)$ is in the category of $\mathcal K/p$--Hopf algebras studied by Bousfield in \cite[Appendix]{bousfield4}.  He shows \cite[Thm.10.8]{bousfield4} that objects in this category are flat over subobjects.  It follows that, if $K(n)_*(\Oinfty f)$ is monic, the spectral sequence associated to the cofibration sequence in $\Alg$
$$\Sinfty (\Oinfty X)_+ \ra \Sinfty (\Oinfty Y)_+
\ra \Sinfty (\Oinfty C)_+$$
collapses, giving the desired short exact sequence.
\end{proof}

\begin{proof}[Proof of \thmref{sn* ses thm}]
Recall that we have cofibration sequences
$$ C_{n-1,d}(X) \xra{f_d} X\langle d \rangle \ra (L_{n-1}^f X)\langle d \rangle,$$
and that $C_{n-1,d}(X) \in \Sp_n^K$ if $d$ is large.

Now consider the diagram used in the proof of \thmref{sn bar iso thm}:
\begin{equation*}
\xymatrix{
K(n)_*(\PP(C_{n-1,d}(X))) \ar[d] \ar[r] &
K(n)_*(\Oinfty C_{n-1,d}(X)) \ar[d]^{(\Oinfty f_d)_*}  \\
K(n)_*(\PP(X \langle d \rangle)) \ar[r]^{s_n(X\langle d \rangle)_*} &
K(n)_*(\Oinfty X \langle d \rangle).
}
\end{equation*}
The left map is always an isomorphism, as is the top map if $d$ is large.  The bottom map is always monic by \thmref{Kn thm}, thus so is the right map, if $d$ is large.  The previous lemma applies to say that, for all large $d$, there is a short exact sequence of $K(n)_*$--Hopf algebras
$$ K(n)_*(\Oinfty C_{n-1,d}(X)) \xra{(\Oinfty f_d)_*} K(n)_*(\Oinfty X\langle d \rangle) \ra K(n)_*(\Oinfty (L_{n-1}^f X)\langle d \rangle).$$
This rewrites as the short exact sequence of the theorem:
$$ K(n)_*(\PP X) \xra{s_n(X)_*} K(n)_*^{vir}(\Oinfty X) \ra K(n)_*^{vir}(\Oinfty L_{n-1}^f X).$$
\end{proof}

\begin{proof}[Proof of \thmref{ses thm}]
Suppose given $f: X\ra Y$ with $X \in \Sp_n^K$.  In the diagram
\begin{equation*}
\xymatrix{
K(n)_*(\PP X) \ar[d]^{\PP(f)_*} \ar[r] &
K(n)_*(\Oinfty X) \ar[d]^{(\Oinfty f)_*}  \\
K(n)_*(\PP Y) \ar[r] &
K(n)_*(\Oinfty Y),
}
\end{equation*}
we then know that the top map is an isomorphism.  Since the bottom map is always monic, if the left map is monic, we deduce that $(\Oinfty f)_*$ is also monic.  If $X$ and $Y$ are 0--connected and $C$ is the cofiber of $f$, the lemma applies, yielding the short exact sequence of the theorem. \\
\end{proof}

\begin{rem} This proof makes evident the following $T(n)_*$ variant of \thmref{ses thm}: if $f: X \ra Y$ is a $T(n)_*$--isomorphism, and $X \in \Sp_n$, then $(\Oinfty f)_*: T(n)_*(\Oinfty X) \ra T(n)_*(\Oinfty Y)$ is monic. \\
\end{rem}

\subsection{Proof of \thmref{sn(X) iso thm} and \thmref{Kn thm 2}}

Since $T(n)_*$ is all $p$--torsion, \cite[\S 4]{bousfield2} implies

\begin{lem} Let $A$ be an abelian group. $K(A,j)$ is $T(n)_*$--acyclic if 

\noindent (1) $j=0$ and $A=0$,

\noindent (2) $1\leq j \leq c(n)-1$ and $A$ is uniquely $p$--divisible, 

\noindent (3) $j=c(n)$ and $A/\text{(torsion)}$ is uniquely $p$--divisible, or

\noindent (4) $j>c(n)+1$.
\end{lem}

\begin{proof}[Proof of \thmref{sn(X) iso thm}(1)]

Given an $S$--module $X$, for each $j \geq 0$, there is a fibration sequence of spaces
$$ K(\pi_{j+1}(X),j) \ra \Oinfty X\langle j+1 \rangle \ra \Oinfty X \langle j \rangle.$$
Under the theorem's hypotheses on $\pi_*(X)$, the fiber is $T(n)_*$--acyclic, so the second map is a $T(n)_*$--isomorphism.  We deduce that $X\langle j+1 \rangle \in \Sp_n$ if and only if $X\langle j \rangle \in \Sp_n$.

The hypothesis that $X \in \bar{\Sp}_n$ means that $X \langle d \rangle \in \Sp_n$ for all large $d$.  By downward induction on $j$, we deduce that $X\langle j \rangle \in \Sp_n$ for all $j \geq 0$.  Since $\pi_0(X) = 0$, $X\langle 0 \rangle \in \Sp_n$ implies $X \in \Sp_n$.

\end{proof}

For $K(n)_*$, we have sharper results.

\begin{lem} Let $A$ be an abelian group. $K(A,j)$ is $K(n)_*$--acyclic if and only if

\noindent (1) $j=0$ and $A=0$,

\noindent (2) $1\leq j \leq n$ and $A$ is uniquely $p$--divisible, 

\noindent (3) $j=n+1$ and $A/\text{(torsion)}$ is uniquely $p$--divisible, or

\noindent (4) $j>n+1$.
\end{lem}

\begin{proof}[Proof of \thmref{sn(X) iso thm}(2) and \thmref{Kn thm 2}]

We can assume that $X \in \bar{\Sp}_n^K$ is 0--connected.  As in our proof of \lemref{strongly acyclic lemma}, let $P^dX$ denote the $d^{th}$ Postnikov section of $X$, so there is a cofibration sequence of spectra
$$ X\langle d \rangle \ra X \ra P^dX.$$
If $d$ is large, then $X \langle d \rangle \in \Sp_n^K$.  Then \thmref{ses thm} applies, and we deduce that there is a short exact sequence of $K(n)_*$--Hopf algebras
$$ K(n)_*(\Oinfty X\langle d \rangle) \ra K(n)_*(\Oinfty X) \ra K(n)_*(\Oinfty P^dX).$$
Thus $X \in \Sp_n^K$ if and only if $\Oinfty P^dX$ is $K(n)_*$--acyclic.

The main theorem of \cite{hrw} says that there is an isomorphism
$$ K(n)_*(\Oinfty P^dX) \simeq \bigotimes_{j=1}^d K(n)_*(K(\pi_j(X),j).$$
\thmref{Kn thm 2} follows.

By the lemma, this tensor product will be isomorphic to $K(n)_*$ if and only if 
$\pi_{j}(X)$ is uniquely $p$--divisible for $1 \leq j \leq n$, and also $\pi_{n+1}(X)/(\text{torsion})$ is uniquely $p$--divisible.  \thmref{sn(X) iso thm}(2) follows.
\end{proof}

\subsection{$s_n(X)$ is universal: proof of \propref{sn is universal prop} and related results}

We prove the first part of \propref{sn is universal prop}; the proofs of the other variants, including \propref{sn* is universal prop}, are similar and left to the reader.

Suppose $F: \Sp \ra \Sp$ is functor preserving $T(n)_*$--isomorphisms, and $T$ is a weak natural transformation of the form
$$ T(X): F(X) \ra L_{T(n)}\Sinfty( \Oinfty X)_+.$$
We show it uniquely factors through $s_n$.

Let $C(X) = C_{n-1,c(n)+2}(X)$ defined as in the proof of \thmref{sn bar iso thm}.  Then $s_n(C(X))$ is an equivalence, and $C(X) \ra X$ is a $T(n)_*$--isomorphism.

We simplify notation; let
$$P(X) = L_{T(n)}\PP(X) \text{ and } L(X) = L_{T(n)}\Sinfty(\Oinfty X)_+.$$
By naturality, we have a commutative diagram
\begin{equation*}
\xymatrix{
F(C(X)) \ar[d]^{\wr} \ar[rr]^{T(C(X))} & &
L(C(X)) \ar[d]  & & P(C(X)) \ar[ll]_{s_n(C(X))}^{\sim} \ar[d]^{\wr}  \\
F(X) \ar[rr]^{T(X)} & &  L(X) & &
P(X) \ar[ll]_{s_n(X)},
}
\end{equation*}
where the left vertical map is an equivalence since $F$ preserves $T(n)_*$--isomorphisms.  The canonical factorization of $T$ through $s_n$ is evident.

\appendix

\section{Comparison of \thmref{classic theorem} with other stable splittings} \label{splitting appendix}

Let $C Z$ be the free $E_{\infty}$--space generated by a space $Z$, as in \cite{may1}.  The inclusion $Z \ra QZ$ then induces the standard approximation map $\alpha(Z): C Z \ra QZ$.  The suspension spectrum $\Sinfty (CZ)_+$ has the structure of an object in  $\Alg$ such that $\Sinfty (\alpha(Z))_+$ is an algebra map.

The purpose of this appendix is to make the following two observations.  Firstly, $s(Z): \PP(\Sinfty Z) \ra \Sinfty (QZ)_+$ refines to a natural map $s^C(X): \PP(\Sinfty Z) \ra \Sinfty (CZ)_+$.  Secondly, $s^C(Z)$ is always an equivalence, and agrees with the standard `stable splittings' in the literature.

This first point is easily checked.  Recall that $s(Z)$ is defined to be the  the natural weak map in $\Alg$ induced by the weak natural map of $S$--modules
$$ \Sinfty Z \xra{\Sinfty \eta(Z)} \Sinfty QZ \xla{\sim} I(\Sinfty(QZ)_+).$$
Similarly we define $s^C(Z)$ to be the natural weak map in $\Alg$ induced by the weak natural map of $S$--modules
$$ \Sinfty Z \xra{\Sinfty \eta(Z)} \Sinfty CZ \xla{\sim} I(\Sinfty(CZ)_+).$$
Then there is an evident factorization
\begin{equation*}
\xymatrix{
\PP(\Sinfty Z) \ar[d]^{s^C(Z)} \ar[dr]^-{s(Z)} &
 \\
\Sinfty(CZ)_+ \ar[r]^-{\alpha(Z)} &
\Sinfty(QZ)_+.
}
\end{equation*}

To check the second point, we begin by observing that $s^C$ admits a slightly different definition.  Let $a(Z)$ denote the fiber (in $S$--modules) of the evident `augmentation' $\Sinfty (Z_+) \ra S$.  Note that the composite $a(Z) \ra \Sinfty (Z_+) \ra \Sinfty Z$ is always an equivalence.  Then $s^C(Z)$ can alternatively be defined as the natural weak map in $\Alg$ induced by the weak natural map of $S$--modules
$$ \Sinfty Z \xla{\sim} a(Z) \ra a(CZ)=I(\Sinfty(CZ)_+).$$

Now we need to recall that $C$ can be defined on the spectrum level \cite{lmms}.  Let $S-\Sp$ denote the category of diagrams of $S$--modules of the form
\begin{equation*}
\xymatrix{
 & X \ar[dr] &  \\
S \ar[ur] \ar@{=}[rr] & &
S.
}
\end{equation*}
There is a functor $C: S-\Sp \ra \Alg$ such that

\noindent (1) $C(\Sinfty(Z_+)) = \Sinfty (CZ)_+$, and 

\noindent (2) $X \ra X \vee S$ induces $\PP(X) = C(X \vee S)$.

Using (2), the commutative diagram
\begin{equation*}
\xymatrix{
\Sinfty Z \ar[d] \ar@{=}[r] &
\Sinfty Z \ar[d] & a(Z) \ar[d] \ar[l]_{\sim} \\
\Sinfty Z \vee S \ar[r]^{\sim} & \Sinfty Z \times S &
\Sinfty (Z_+) \ar[l]_{\sim},
}
\end{equation*}
induces a diagram in $\Alg$
\begin{equation*}
\xymatrix{
\PP(\Sinfty Z) \ar@{=}[d] \ar@{=}[r] &
\PP(\Sinfty Z) \ar[d] & \PP(a(Z)) \ar[d] \ar[l]_{\sim} \\
C(\Sinfty Z \vee S) \ar[r]^{\sim} & C(\Sinfty Z \times S) &
C(\Sinfty (Z_+)) \ar[l]_{\sim}.
}
\end{equation*}
Now using (1), this shows that $s^C(Z)$ is the natural weak equivalence $$ \PP(\Sinfty Z) = C(\Sinfty Z \vee S) \xra{\sim} C(\Sinfty Z \times S) \xla{\sim} \Sinfty (CZ)_+.$$
Defined this way, $s^C(Z)$ satisfies the characterization of natural splittings given in \cite[Appendix B]{k3}. \\

We end this appendix by noting that the proof of \thmref{susp thm} generalizes in a straightforward way to prove the following variant. \\

\begin{thm} \label{Csusp thm}  If a map of spectra $f: \Sinfty Z \ra X$ is an $E_*$--isomorphism, then the composite
$$ E_*(CZ) \xra{\alpha(Z)} E_*(Q Z) \xra{(\Oinfty f)_*} E_*(\Oinfty X)$$
is a monomorphism. \\
\end{thm}

\section{Comparison with recent work of Bousfield} \label{Bousfield appendix}

In this appendix, we show how \thmref{classic theorem} and  \thmref{Tn theorem} can be compared by  using Bousfield's beautiful natural zig--zag of $L_{T(n)}$--equivalences relating any $S$--module $X$ to a suspension spectrum determined by $X$ \cite{bousfield6}.  This allows for an alternative proof of \thmref{sn(X) iso thm 2}, and thus of many of the results in \secref{sn iso subsection} and \secref{more homo cor subsection}.   

Bousfield constructs a functor
$$ \Theta_n: \text{$S$--modules} \ra \text{Spaces}$$
that is a left adjoint of sorts to the telescopic functor 
$$ \Phi_n: \text{Spaces} \ra \text{$S$--modules}.$$
Using this adjunction, the equivalence 
$$L_{T(n)}X \xra{\sim} \Phi_n(\Oinfty X),$$
 corresponds to an equivalence
$$L_n^f \Sinfty \Theta_n(X) \xra{\sim} M_n^fX,$$ 
where $M_n^f X$ is the fiber of $L_n^f X \ra L_{n-1}^f X$.  Thus $L_{n-1}^f \Sinfty \Theta_n(X) \simeq *$, and there is a natural $T(n)_*$--equivalence $\Sinfty \Theta_n(X) \ra L_{n}^f X$. Furthermore, Bousfield observes that $\Theta_n(X)$ is always $d_n$--connected, where $d_n$ is defined in \cite[\S 4.3]{bousfield6}: one can deduce that $d_n \geq c(n)+1$ from \cite[Prop.2.1]{bousfield2}.  One also has \cite[Thm.3.3]{bousfield6} that $M_n^f L_{T(n)} \simeq M_n^f$ and $L_{T(n)} M_n^f \simeq L_{T(n)}$. 

Thus we get a zig--zag of $T(n)_*$--equivalences:
$$ \Sinfty \Theta_n(X) \xra{\beta(X)} (L_{n}^f X)\langle d_n \rangle \la X\langle d_n \rangle  \ra X.$$

This allows us to consider the following diagram:

\begin{equation*}
\xymatrix{
L_{T(n)}\PP(X) \ar[rrr]^-{s_n(X)} & & & L_{T(n)}\Sinfty \Oinfty X_+   \\
L_{T(n)}\PP(X\langle d_n \rangle) \ar[rrr]^-{s_n(X\langle d_n \rangle)} \ar[d]_{\wr}  \ar[u]^{\wr} && & L_{T(n)}\Sinfty \Oinfty X\langle d_n \rangle_+ \ar[d]_{\wr} \ar[u] \\
L_{T(n)}\PP((L_n^fX)\langle d_n \rangle)  \ar[rrr]^-{s_n((L_n^fX)\langle d_n \rangle)} & && L_{T(n)}\Sinfty \Oinfty (L^f_nX)\langle d_n \rangle_+   \\
L_{T(n)}\PP(\Sinfty \Theta_n(X)) \ar[u]^{\wr} \ar[rrr]^-{L_{T(n)}s(\Theta_n(X))}_{\sim} & &&  L_{T(n)}\Sinfty Q \Theta_n(X)_+. \ar[u]  \\
}
\end{equation*}

Below we will show that the diagram commutes.  Thus the classical stable splitting of $Q \Theta_n(X)$ given by \thmref{classic theorem}, corresponds to the splitting of $L_{T(n)}\Sinfty \Oinfty X$ given by \thmref{Tn theorem}.

A crucial point about this diagram is that, as indicated, the middle vertical map on the right is an equivalence, as $\Oinfty$ takes $L_n^f$--equivalences between $d_n$--connected spectra to $T(n)_*$--equivalences, thanks to \cite[Cor.4.8]{bousfield6}\footnote{As $d_n \geq c(n)+1$, this also follows from \thmref{sn(X) iso thm 2}. However, if one wishes to offer an alternative proof of \thmref{sn(X) iso thm 2}, it seems prudent to not argue this way.}.   Thus, since the diagram commutes, it is clear that $s_n(X)$ is an equivalence on highly connected $X$ if and only if the bottom right vertical map is an equivalence, i.e. $(\Oinfty \beta(X))_*: T(n)_*(Q \Theta_n(X)) \ra T(n)_*(\Oinfty (L_{n}^f X)\langle d_n \rangle)$ is an isomorphism.  Again appealing  to \cite[Cor.4.8]{bousfield6}, this will happen if $L_{n-1}^f X \simeq *$, and we have reproved \thmref{sn(X) iso thm 2} using Bousfield's results.

It is illuminating to note that $\Oinfty \beta(X)$ is always a $T(n)_*$--monomorphism, by virtue of our \thmref{susp thm}.

The top two squares of the diagram obviously commute.  Checking that the bottom square commutes quickly reduces to verifying that the following diagram commutes:
\begin{equation*}
\xymatrix{
L_{T(n)}(L_n^fX)\langle d_n \rangle  \ar[rrr]^-{\eta_n((L_n^fX)\langle d_n \rangle)} & && L_{T(n)}\Sinfty \Oinfty ((L_n^fX)\langle d_n \rangle)   \\
L_{T(n)}\Sinfty \Theta_n(X) \ar[u]^{\wr}_{L_{T(n)}\beta(X)} \ar[rrr]^-{L_{T(n)}\Sinfty \eta(\Theta_n(X))} & &&  L_{T(n)}\Sinfty Q \Theta_n(X). \ar[u]_{L_{T(n)}\Sinfty \Oinfty \beta(X)}  \\
}
\end{equation*}

We show this using a variant of a proof which was outlined to us in email from Pete Bousfield.  It is an exercise in using the various adjunctions constructed in \cite{bousfield6}, as summarized in \cite[Thm.5.14]{bousfield6}.

By the naturality of $\eta_n$, it suffices to verify the following proposition. \\

\begin{prop} \label{theta prop}  $ \eta_n(\Sinfty \Theta_n(X)) \simeq L_{T(n)}\Sinfty \eta(\Theta_n(X))$. \\
\end{prop}

To prove this, we first observe that $\Sinfty \Theta_n$ preserves $T(n)_*$--equivalences, and thus so does $\Sinfty Q \Theta_n$.  Thus the zig--zag of $T(n)_*$--equivalences
$$ X \ra L_n^f X \la M_n^f X$$
can be used to reduce the proof of the proposition to the case when $X = M_n^f X$, i.e. $X \in \mathcal M_n^f$ in the notation of \cite{bousfield6}.

For any space $Z$,  unravelling the definitions reveals that 
$$\eta_n(\Sinfty Z) = \Phi_n(\eta(QZ)), $$
while 
$$ L_{T(n)} \Sinfty \eta(Z) = \Phi_n(Q\eta(Z)).$$
Both of these maps clearly agree after precomposition with 
$$ \Phi_n(\eta(Z)): \Phi_n(Z) \ra L_{T(n)} \Sinfty Z.$$
Thus the next lemma will compete the proof of the proposition. \\

\begin{lem} \label{theta lemma}  If $X \in \mathcal M_n^f$, then 
$$ \Phi_n(\eta(\Theta_n(X))): \Phi_n(\Theta_n(X)) \ra L_{T(n)} \Sinfty \Theta_n(X)$$
is split epic.
\end{lem}

To prove this, we will use a very general categorical lemma.  Suppose one has two categories $\A$ and $\B$, and two pairs of adjoint functors
\begin{equation*}
\xymatrix{
\A  \ar@<1ex>[r]^{L_1} &  \B \ar@<1ex>[r]^{L_2} \ar@<1ex>[l]^{R_1} & \A  \ar@<1ex>[l]^{R_2}.
}
\end{equation*}
Let $\eta_1: 1_{\A} \ra R_1 L_1$ and $\eta_2: 1_{\B} \ra R_2L_2$ be the units of the adjunctions.

Now suppose we are also given a natural transformation $\gamma: 1_{\A} \ra R_1R_2$ with adjoint $\beta: L_2 L_1 \ra 1_{\A}$.  

\begin{lem} $\gamma$ is an equivalence if and only if $\beta$ is an equivalence. In this case, the map
$$ R_1 \eta_2(L_1A): R_1 L_1 A \ra R_1 R_2 L_2 L_1 A$$
is split epic for all $A$.
\end{lem}

\begin{proof}  The first statement is clear.  The second statement then follows from the commutative diagram
\begin{equation*}
\xymatrix{
A \ar[d]^{\gamma(A)}_{\wr} \ar[rr]^-{\eta_1(A)} &&
R_1L_1A \ar[d]^{R_1\eta_2(L_1A)}  \\
R_1R_2A  &&  R_1R_2L_2L_1A \ar[ll]^{\sim}_{R_1R_2\beta(A)}.
}
\end{equation*}
\end{proof}

\begin{proof}[Proof of \lemref{theta lemma}]  The previous lemma applies to the pair of adjoint functors appearing in \cite[Thm.5.14]{bousfield6} to say that
$$ M_n^f \Phi_n(\eta_2(\Theta_n(X))): M_n^f \Phi_n\Theta_n(X) \ra M_n^f \Phi_n \Oinfty (L_n^f \Sinfty \Theta_n(X))\langle d_n \rangle$$
is split epic, where $\eta_2$, defined on a certain category of $d_n$--connected spaces, has the form
$$ \eta_2(Z): Z \ra \Oinfty(L_n^f \Sinfty Z)\langle d_n \rangle.$$

Applying $L_{T(n)}$, one deduces that 
$$ \Phi_n(\eta_2(\Theta_n(X))): \Phi_n \Theta_n(X) \ra L_{T(n)}((L_n^f \Sinfty \Theta_n(X))\langle d_n \rangle)$$
is split epic.

Using the zig--zag of $T(n)_*$--equivalences
$$ (L_n^f \Sinfty \Theta_n(X))\langle d_n \rangle \ra L_n^f \Sinfty \Theta_n(X) \la \Sinfty \Theta_n(X),$$
it follows that this last map identifies with
$$ \Phi_n(\eta(\Theta_n(X))): \Phi_n \Theta_n(X) \ra L_{T(n)} \Sinfty \Theta_n(X).$$
\end{proof}

\section{A comparison of $s_n$ and $s$} \label{hopf invariant appendix}

One might wonder for what $Z$ the two natural maps
$$s_n(\Sinfty Z), L_{T(n)}s(Z): L_{T(n)} \PP(\Sinfty Z) \ra L_{T(n)} \Sinfty (QZ)_+$$
are homotopic.  Here we briefly summarize what we can say about this.

On the positive side, $s_n(\Sinfty Z) \simeq L_{T(n)}s(Z)$ if and only if 
$$\eta_n(\Sinfty Z) \simeq L_{T(n)}\eta(Z): L_{T(n)} \Sinfty Z \ra L_{T(n)} \Sinfty QZ.$$
Thus \propref{theta prop} implies \\

\begin{prop} \label{theta prop 2} If $Z \simeq \Theta_n(X)$ then $s_n(\Sinfty Z) \simeq L_{T(n)}s(Z)$. \\
\end{prop}

\lemref{s = sn lem} gave another sufficient condition on $Z$; we do not know if the proposition includes this as a special case.

On the negative side, since $s(Z)$ is an equivalence for all connected $Z$, we have an obvious necessary condition. \\

\begin{lem} If $Z$ is connected and $\Sinfty Z \not\in \Sp_n$ then $s_n(\Sinfty Z) \not\simeq L_{T(n)}s(Z)$.
\end{lem}

Thus, for example, the two maps are distinct for $Z = S^1$, and for all $Z=S^d$ if $n\geq 2$.   More examples of suspension spectra not in $\Sp_n$ can be found using the next simple lemma, which doesn't seem to follow immediately from our  other results. \\

\begin{lem} If $X \not\in \Sp_n$ then $\Sinfty \Oinfty X \not\in \Sp_n$.
\end{lem}

\begin{proof}  The evaluation map $\Sinfty \Oinfty X \ra X$ induces a commutative diagram
\begin{equation*}
\xymatrix{
T(n)_*(\PP \Sinfty \Oinfty X) \ar[d] \ar[rr]^{s_n(\Sinfty \Oinfty X)} &&
T(n)_*(Q \Oinfty X) \ar[d]  \\
T(n)_*(\PP X) \ar[rr]^{s_n(X)} &&
T(n)_*(\Oinfty X).
}
\end{equation*}
The horizontal maps are always monic, and the right vertical map is always epic, as it admits an obvious splitting.  Thus if the top map is an isomorphism, so is the bottom.
\end{proof}

When $\eta_n(\Sinfty Z)$ and $L_{T(n)}\eta(Z)$ differ, one can roughly measure the difference by means of James--Hopf invariants.  Let
$$ t_r(Z): \Sinfty Q Z \ra \Sinfty D_r Z$$
be the $r^{th}$ component of $t(Z)$, as given by \thmref{classic theorem}.
In the literature, the adjoint
$$ j_r(Z): QZ \ra QD_rZ$$
is usually called the $r^{th}$ James--Hopf invariant.

Now consider the two maps
$$ L_{T(n)} \Sinfty Z \xra{L_{T(n)}\eta(Z)} L_{T(n)} \Sinfty QZ \xra{L_{T(n)}t_r(Z)} L_{T(n)} \Sinfty D_r Z,$$
and 
$$ L_{T(n)} \Sinfty Z \xra{\eta_n(\Sinfty Z)} L_{T(n)} \Sinfty QZ \xra{L_{T(n)}t_r(Z)} L_{T(n)} \Sinfty D_r Z.$$

For $r \geq 2$, the former is 0, while the latter is $\Phi_n(j_r(Z))$, as is easily checked.

Comparison with \propref{theta prop 2} implies the next corollary. 

\begin{cor} If $s_n(\Sinfty Z) \simeq L_{T(n)}s(Z)$, e.g. if $Z \simeq \Theta_n(X)$ for some $X$, then $\Phi_n(j_r(Z))$ is null for all $r \geq 2$.
\end{cor}

There are some intriguing open questions regarding the natural transformations
$$ \Phi_n(j_r(Z)): L_{T(n)} \Sinfty Z \ra L_{T(n)} \Sinfty D_r Z.$$

For example, they induces natural transformations
$$ E_n^*(D_r Z) \ra E_n^*(Z),$$
and one might wonder if these are related to either the constructions in \cite{hkr}, or, via duality, to total power operations in $E_n^*$ cohomology.

\end{document}